\documentclass{article}

\usepackage{times}
\usepackage{amssymb}

\providecommand{\keywords}[1]
{
  \small	
  \textbf{\textit{Keywords---}} #1
}

\title{Meaning as Use, Application, Employment, Purpose, Usefulness\thanks{Preprint. Final version to be published in De Gruyter's \emph{SATS: Northern European Journal of Philosophy.}}}

\author{
Ruy J.G.B. de Queiroz\\
Centro de Inform\'atica\\
Universidade Federal de Pernambuco\\
Recife, PE, Brazil\\
ruy@cin.ufpe.br}

\begin{document}

\maketitle

\begin{abstract}
Arising from the whole body of Wittgenstein's writings is a picture of a (not necessarily straight, linear, but admittedly tireless) journey to come to terms with the mechanics of language as an instrument to conceive `reality' and to communicate an acquired conception of the `world'. The journey passes through mathematics, psychology, color perception, certainty, aesthetic, but, looking at it from a sort of birdview, it seems reasonable to say that these are all used as `test beds' for his reflections and `experimentations' towards an all encompassing perspective of such a fundamental gateway to human reasoning and life-revealing as language. Whatever labelling of Wittgenstein as a mystic, a logicist, a conventionalist, a skeptic, an anti-metaphysics, an anti-realist, a verificationist, a pragmatist, a behaviorist, and many others, does not seem to do justice to his absolute obsession with being a persistent `deep diver' into the nature of language.
Working with an open and searchable account of the \emph{Nachlass} has allowed us to identify important aspects of the philosopher's possible common line of thinking, in spite of changes of directions, some of them acknowledged by Wittgenstein himself. One of those aspects is the association of meaning with use, application, purpose, usefulness of symbols in language, which happens to show itself from the very beginning through to the very late writings. The German terms \emph{Gebrauch}, \emph{Anwendung}, \emph{Verwendung}, \emph{Zweck} in relation to meaning, sense of signs, words, sentences, appear in several texts since the \emph{WW1 Notebooks (1914--1916)} up until very late manuscripts from 1950--51.\footnote{This is a followup to two recent papers, one already published in 2023, and another one which has just been accepted for publication. We have already made heavy use of Wittgenstein's \emph{Nachlass}, and continued to do so here, much more than it was done in any of the previous publications. Thus, as a new step in a series of essays going back a few decades, it certainly contains some overlap with previous publications intended to bring context and self-containedness to the present manuscript, this substantially reinforces earlier arguments. (Thanks to the editors of \emph{Wittgenstein Archives at the University of Bergen} (\emph{WAB}) who have awarded Wittgenstein scholarship with such a huge gift as the \emph{Nachlass} in free online and searchable form!)}

\end{abstract}
 

\section{Introduction}
The use of the open and searchable Wittgenstein's \emph{Nachlass} (The Wittgenstein Archives at the University of Bergen (WAB))  has served spectacularly in the quest for what seems to constitute a common thread of Wittgenstein's view on the connections between meaning, use, purpose, and consequences, going from the \emph{Notebooks} to later writings (including and beyond the \emph{Philosophical Investigations}). And this search for a baseline of Wittgenstein's thinking on language looks for taking it as the basis for a proposal for a formal counterpart of a `meaning-as-use' (dialogical/game-theoretical) semantics for the language of predicate logic. In order to further consolidate a perspective first raised back four decades ago, we have been bringing out key excerpts from Wittgenstein's \emph{Nachlass}, especially from the impressive repository of WAB.
Having Wittgenstein's \emph{Nachlass} at hand with the possibility of working with computer-assisted search has proven extremely useful in reexamining not just the unpublished writings, but also the body of published work, to the extent that it has allowed us to identify important aspects of the philosopher's possible common line of thinking, in spite of changes of directions, some of them acknowledged by Wittgenstein himself. One of those aspects is the association of meaning with use, application, purpose, usefulness of symbols in language, which happens to show itself from the very beginning through to the very late writings. The German terms \emph{Gebrauch}, \emph{Anwendung}, \emph{Verwendung}, \emph{Zweck} in relation to meaning, sense of signs, words, sentences, appear in several texts since the \emph{WW1 Notebooks (1914--1916)} up until very late manuscripts from 1950--51, including (and beyond) the published text of the \emph{Investigations}. In this context, it seems appropriate to refer to what the Editorial Preface of the revised fourth edition of the \emph{Investigations},
by P. M. S. Hacker and Joachim Schulte says:
\begin{quote}
Anscombe was not consistent in her translation of \emph{Gebrauch}, \emph{Verwendung} and \emph{Anwendung}. We have translated \emph{Gebrauch} by `use', \emph{Verwendung} by `use' or `employment', and \emph{Anwendung} by `application'. `Use' also does service for \emph{ben\"utzen}. In general, however, we have not allowed ourselves to be hidebound by the multiple occurrence of the same German word or phrase in different contexts. It by no means
requires always translating by the same English expression, but rather
depends on the exigencies of the context and the author's intention. So,
for example, we have translated \emph{Praxis der Sprache} in \emph{Investigations}
\S\ 21 by `linguistic practice' rather than by the more ponderous `practice of the language', and \emph{Praxis des Spiels} in \S\ 54(b) as `the way the
game is played', because this is how Wittgenstein wanted it translated.
\end{quote}

By examining those passages quoted here from his unpublished as well as published works, we wish to contend that the results of the search in the \emph{Nachlass} strongly support and augment the main thrust of our continuing research on the significance of Wittgenstein's suggestion that meaning is revealed by explaining the purpose, the usefulness, the immediate consequences of a term or statement. In particular, we have recently shown in two previous papers (2023, published, and 2024, under review) that in spite of the dominant view of the picture theory of the \emph{Tractatus}, one finds clear signs in very early (pre-Tractatus) writings of what became the so-called `meaning as use' paradigm of the second phase. Some of those signs had not been explicitly uncovered before, and this was made possible due to the openness and the easiness of making computer-assisted search in the whole body of writings.
By examining these new findings in the \emph{Nachlass}, we have sought to connect Wittgenstein's account of 'meaning as use', specified by the 'calculus' involved in specifying or using terms (or statements), with the pragmatist approaches to meaning as revealed in Peirce's writings (1932) on the interaction between the \emph{Interpreter} and the \emph{Utterer} which seem to bear on the ``game"/``dialogue" approaches to meaning (Lorenzen, Hintikka, Fra\"{\i}ss\'e). And this seems to be compatible with Wittgenstein's own attitude towards the possible labelling him as a pragmatist, which one finds in a late manuscript (Ms-131, 1946) but not in the published \emph{Investigations}:
\begin{quote}
Wie aber, wenn die Religion lehrt, die Seele k\"onne bestehen, wenn der Leib zerfallen ist? Verstehe ich, was sie lehrt? Freilich verstehe ich's: Ich || ich kann mir dabei manches vorstellen.

(Man hat ja auch Bilder von diesen Dingen gemalt. Und warum sollte so ein Bild nur die unvollkommene Wiedergabe des ausgesprochenen Gedankens sein? Warum soll es nicht den gleichen Dienst tun, wie der Satz? -- wie das, was wir sagen?)

     Und auf den Dienst kommt es an.
 	 			 	
     Aber bist Du kein Pragmatiker? Nein. Denn ich sage nicht, der Satz sei wahr, der n\"utzlich ist.
     
     Der Nutzen, d.h., Gebrauch, gibt dem Satz seinen besondern Sinn, das Sprachspiel gibt ihm ihn.
     
     Und insofern als eine Regel oft so gegeben wird, da\ss\ sie sich n\"utzlich erweist, \& mathematische S\"atze mit Regeln || ihrer Natur || ihrem Wesen nach mit Regeln verwandt sind, spiegelt sich in mathematischen
Wahrheiten N\"utzlichkeit.\footnote{which can be translated as:
\begin{quote}
But what if religion teaches that the soul can exist when the body has decayed? Do I understand what she teaches? Of course I understand it: I || I can imagine a lot of things.

(People have also painted pictures of these things. And why should such a picture only be an imperfect representation of the thought expressed? Why shouldn't it do the same job as the sentence? || , as what we say?)
     
     And the service is what matters.

      But aren't you a pragmatist? No. For I do not say that the sentence that is useful is true.
      
      Utility, i.e. use, gives the sentence its special meaning; the language game gives it this.
      
      And insofar as a rule is often given in such a way that it proves useful, \& mathematical sentences with rules || their nature || are essentially related to rules, usefulness is reflected in mathematical truths.
\end{quote}
Ms-131,69-71 (1946),  Ts-245,185 (1945--47), and Ts 229,253 (1947) Wittgenstein Nachlass Ms-131 (WL).
User filtered transcription. In: Wittgenstein, Ludwig: Interactive Dynamic Presentation (IDP) of Ludwig Wittgenstein's philosophical Nachlass. Edited by the Wittgenstein Archives at the University of Bergen (WAB) under the direction of Alois Pichler. Bergen: Wittgenstein Archives at the University of Bergen 2016-. (Accessed 02 Jan 2025)
}
\end{quote}

And, indeed, the unrestricted access to a searchable \emph{Nachlass} may allow for a reevaluation of certain widespread impressions, in most cases without access to the body of the \emph{Nachlass}, that Wittgenstein, especially in the later period, should be seen as a conventionalist with respect to logical necessity, as, e.g., Dummett argued on his discussions on the debate regarding realism vs anti-realism. By taking into account the wider spectrum of writings on language understanding and conception of reality, which arises out of a deep dive in the posthumous (unpublished) writings, one can see that Wittgenstein was consistently focussed on seeking to unravel the process of making the connections between language and `reality', though bumping into self-corrections, even if the actual nature of those connections may or may not be taken to be realist or anti-realist as this dichotomy underlies some discussions around language and the nature of mathematics. In particular,
\begin{quote}
In der Logik geschieht immer wieder, was in dem Streit \"uber das Wesen der Definition geschehen ist. Wenn man sagt, die Definition habe es nur mit Zeichen zu tun und ersetze blo\ss\ ein kompliziertes Zeichen durch ein einfacheres -- ein Zeichen durch ein anderes, so wehren sich die Menschen dagegen und sagen, die Definition leiste nicht nur das, oder es gebe eben verschiedene Arten von Definitionen -- der Definition und die interessante und wichtige sei nicht die (reine) ``Verbaldefinition".

     Sie glauben n\"amlich, man nehme der Definition ihre Bedeutung, Wichtigkeit, wenn man sie als blo{\ss}e Ersetzungsregel, die von Zeichen handelt, hinstellt. W\"ahrend die Bedeutung der Definition in ihrer Anwendung liegt, quasi in ihrer Lebenswichtigkeit. Und eben das geht (heute) in dem Streit zwischen Formalismus, Intuitionismus, etc. vor sich. Es ist den Leuten? unm\"oglich, die Wichtigkeit einer Sache || Handlung || Tatsache, ihre Konsequenzen, ihre Anwendung, von ihr selbst zu unterscheiden; die Beschreibung einer Sache von der Beschreibung ihrer Wichtigkeit.\footnote{which can be translated as:
\begin{quote}
In logic, what happened in the argument about the nature of the definition happens again and again. When it is said that the definition only deals with signs and merely replaces a complicated sign with a simpler one, people object and say that the definition does not only do that, or that there are different kinds of definitions of the definition and that the interesting and important one is not the (pure) `verbal definition'.

     They believe that the definition is deprived of its meaning and importance if it is presented as a mere substitution rule that deals with signs. Whereas the meaning of the definition lies in its application, in its vitality, so to speak. And this is precisely what is going on (today) in the dispute between formalism, intuitionism, etc. It is impossible for people? to distinguish the importance of a thing || action || fact, its consequences, its application, from itself; the description of a thing from the description of its importance.
\end{quote}
Ms-112,18r (1931) (Ms-112: VIII, Bemerkungen zur philosophischen Grammatik (\"ONB, Cod.Ser.n.22.021)) (Accessed 02 Jan 2025)
}
\end{quote}
In a manuscript written in 1931 (Ms-111), there is a direct reference to what seems to be a central point of debate between realism and anti-realism, i.e.\ proof of existence, as in:
\begin{quote}
 ,,Jeder Existenzbeweis mu\ss\ eine Konstruktion dessen enthalten dessen Existenz er beweist." Man kann nur sagen ,,ich nenne Existenzbeweis nur einen der eine solche Konstruktion enth\"alt". Der Fehler ist || liegt darin da\ss\ man glaubt || vorgibt einen klaren Begriff des Existenzbeweises || der Existenz zu besitzen.
 
     Man glaubt ein Etwas, die Existenz, beweisen zu k\"onnen, soda\ss\ man nun unabh\"angig vom Beweis von ihr \"uberzeugt ist. (Die Idee der voneinander ? und daher wohl auch vom Bewiesenen ? unabh\"angigen Beweise!) In Wirklichkeit ist Existenz das was man mit dem beweist, was man ,,Existenzbeweis? nennt. Wenn die Intuitionisten \& Andere dar\"uber reden so sagen sie: ,,Diesen Sachverhalt, die Existenz, kann man nur so \& nicht so beweisen." Und sehen nicht, da\ss\ sie damit einfach das definiert haben was sie Existenz nennen. Denn die Sache verh\"alt sich eben nicht so wie wenn man sagt: ,,da\ss\ er || ein Mann in dem Zimmer ist kann man nur dadurch beweisen da\ss\ man hineinschaut aber nicht, indem man an der T\"ure horcht"
 	 			 	
Wir haben keinen Begriff der Existenz unabh\"angig von unserem Begriff des Existenzbeweises.
 	 			 	
17 + 28 kann ich mir nach Regeln ausrechnen ich brauche 17 + 28 = 45($\alpha$) nicht als Regel zu geben. Kommt also in einem Beweis der \"Ubergang von f(17 + 28) auf f(45) vor so brauche ich nicht sagen er gesch\"ahe nach der Regel $\alpha$, sondern nach andern Regeln des 1 + 1.\footnote{which can be translated as:
\begin{quote}
`Every proof of existence must contain a construction of whose existence it proves.' One can only say `I call proof of existence only one that contains such a construction'. The mistake is that one believes || pretends to have a clear concept of the proof of existence || of existence.

     One believes to be able to prove something, existence, so that one is now convinced of it independently of the proof. (The idea of proofs that are independent of each other - and therefore also of what is proved!) In reality, existence is what one proves with what is called ?proof of existence?. When the intuitionists \& others talk about it, they say: ?This fact, existence, can only be proven in this way \& not in that way.? And don't realise that they have simply defined what they call existence. For the matter does not behave in the same way as when one says: ?That he || a man is in the room can only be proved by looking in but not by listening at the door?.
 	 			 	
We have no concept of existence independent of our concept of proof of existence.

7 + 28 can be calculated according to rules I do not need to give 17 + 28 = 45($\alpha$) as a rule. So if the transition from f(17 + 28) to f(45) occurs in a proof, I do not need to say it happens according to the rule $\alpha$, but according to other rules of 1 + 1.
\end{quote}
Ms-111,155 (Ms-111: VII, Bemerkungen zur Philosophie (WL)) (Accessed 26 Dec 2024)
}
\end{quote}

In spite of a very early conception of meaning as related to purpose, use, usefulness, as we will see in the sequel, the influence of logic via Russell, Frege and other references in the tradition, led Wittgenstein to a search for what he himself referred to as `phenomenological language' which resulted in the so-called picture theory with the notion of `state of affairs' as the basis for it. As in several occasions in his later writings, there comes a self-correction in a form which we would like to call a `go back to his roots' in the sense of meaning as use, purpose, usefulness, application:
\begin{quote}
``Ph\"anomenologische Sprache." Glaube an ihre Notwendigkeit. Es schien als sei unsere Sprache, irgendwie, roh, eine unvollkommene Darstellung der Sachverhalte \& nur als rohes, unvollkommenes Abbild zu verstehen. Als m\"u{\ss}te die Philosophie sie verbessern, verfeinern, um so den Bau der Welt verstehen zu k\"onnen. Dann wurde es offenbar da\ss\ sie die Sprache wie sie ist verstehen, d.h. erkennen m\"usse, weil nicht eine neue Klarheit, die die alte Sprache nicht gebe || gibt, das Ziel sei, sondern die Beseitigung der philosophischen { Irrg\"arten, bewilderment.

     Perplexities, Herumirren, Ratlosigkeit, Sich-nicht-auskennen, Perplex-sein, R\"atsel, Irrfragen, ich meine eine Frage die dazu gestellt ist in einen Irrgarten zu f\"uhren, wie sie in etwa eine Hexe stellen k\"onnte. }
     
     ``Und Deine Schwierigkeiten sind Mi{\ss}verst\"andnisse." ? Wenn sie nicht Mi{\ss}verst\"andnisse w\"aren, wenn wir wirklich weiter analysieren m\"u{\ss}ten um auf festen || sichern Grund \& Boden zu kommen, dann m\"u{\ss}ten wir uns fragen: durch welchen feineren Begriff haben wir den landl\"aufigen Begriff ``Wort" (z.B.) zu ersetzen. Wir m\"u{\ss}ten dann die \"ublichen W\"orter || Ausdr\"ucke mit ihrem Gebrauch zur Seite r\"aumen um in der Tiefe nach den eigentlichen Begriffen zu forschen || suchen nach denen wir die Sprache anpassen m\"ussen || anzupassen || einzurichten haben || h\"atten.\footnote{which can be translated as:
\begin{quote}
`Phenomenological language.' Belief in its necessity. It seemed as if our language was somehow, crudely, an imperfect representation of the facts \& could only be understood as a crude, imperfect image. As if philosophy had to improve it, refine it, in order to be able to understand the structure of the world. Then it became obvious that it had to understand the language as it is, i.e. recognise it, because the goal was not a new clarity, which the old language does not give, but the elimination of the philosophical {misguidedness, bewilderment.

     Perplexities, wandering, perplexity, not knowing, being perplexed, riddles, wrong questions, I mean a question that is asked to lead into a maze, such as a witch might ask. }
     
     `And your difficulties are misunderstandings.' - If they were not misunderstandings, if we really had to analyse further in order to arrive at solid ground, then we would have to ask ourselves: with what finer term do we have to replace the common term `word' (for example)? We would then have to set aside the usual words || expressions and their usage in order to search in depth for the actual terms according to which we would have to adapt the language.
\end{quote}
Ms-152,92 (1936) (Wittgenstein Nachlass Ms-152 [so-called C8] (WL) (Accessed 26 Dec 2024)
}
\end{quote}
Going back to 1914--1915 (with the \emph{Notebooks}), one finds a clear manifestation of what is considered as part of the later paradigm (meaning, use, consequences, purpose, application):
\begin{quote}
23.10.14\\
Um das Zeichen im Zeichen zu erkennen mu\ss\ man auf den Gebrauch achten.\footnote{which can be translated as:
\begin{quote}
23.10.14\\
To recognize the sign in the sign, you have to pay attention to the use.
\end{quote}
Ms-101,59r (Ms-101 [so-called WW1 notebooks] (WL)) (Accessed 02 Jan 2025)
}

24.11.14\\
Da\ss\ man aus dem Satz ,,$(x) \phi x$" auf den Satz ,,$\phi a$" schlie{\ss}en kann das zeigt wie die Allgemeinheit auch im Zeichen ,,$(x) \phi x$" vorhanden ist.
Und das gleiche gilt nat\"urlich f\"ur die Allgemeinheitsbezeichnung \"uberhaupt.\footnote{which can be translated as:
\begin{quote}
24.11.14\\
That the proposition ``$\phi a$" can be inferred from the proposition ``$(x) \phi x$" shews how generality is present even in the sign ``$(x) \phi x$". And the same applies, of course, to the notation of generality in general.
\end{quote}
Ms-102,37r (Ms-102 [so-called WW1 notebooks] (WL)) (Accessed 02 Jan 2025)
}

23.4.15\\
Es ist klar da\ss\ Zeichen, die denselben Zweck erf\"ullen logisch identisch sind. Das rein Logische ist eben das was alle diese leisten k\"onnen.\footnote{which can be translated as:
\begin{quote}
23.4.15
It is clear that signs that fulfil the same purpose are logically identical. The purely logical is precisely what they can all do.
\end{quote}
Ms-102,75r (Ms-102 [so-called WW1 notebooks] (WL)) (Accessed 02 Jan 2025)
}

30.5.15\\
Namen kennzeichnen die Gemeinsamkeit einer Form und eines Inhalts. - Sie kennzeichnen erst mit ihrer syntaktischen Verwendung zusammen eine bestimmte logische Form.\footnote{which can be translated as:
\begin{quote}
30.5.15\\
Names characterise the commonality of a form and a content. - They only characterise a certain logical form together with their syntactic use.
\end{quote}
Ms-102,116r (Ms-102 [so-called WW1 notebooks] (WL)) (Accessed 02 Jan 2025)
}
\end{quote}
The association of meaning to use, purpose, application, via Gebrauch, Zweck, Anwendung, Verwendung, is retaken with full force in the early second phase, as Ms-107 from 1929--30 shows:
\begin{quote}
Wer es gesehen h\"atte, h\"atte gesehen wie kompliziert es ist \& seine Komplikation erkl\"art sich nur durch den beabsichtigten Gebrauch zu dem es tats\"achlich nicht gekommen ist. So m\"ochte ich bei der Sprache sagen: Wozu alle diese Ans\"atze, sie haben nur dann eine Bedeutung wenn sie Verwendung finden.
			 	
Kann man sagen: der Sinn des || eines Satzes ist sein Zweck? [Oder von einem Wort ,,its meaning is its purpose".]
			 	
Die Logik kann aber nicht die Naturgeschichte des Gebrauchs eines Worts angehen.
			 	
In wiefern kann der Gebrauch einer Form eines Wortes die Existenz der || einer anderen Form voraussetzen? So w\"urde etwa ,,dunkel" ,,dunkler" voraussetzen und umgekehrt; oder ,,wei{\ss}" ,,wei{\ss}lich" und umgekehrt etc.\footnote{which can be translated as:
\begin{quote}
Anyone who had seen it would have seen how complicated it is \& its complication can only be explained by the intended use to which it was not actually put. So I would like to say with language: why all these approaches, they only have a meaning when they are used.
 	 			 	
Can one say: the meaning of a sentence is its purpose? [Or of a word ?its meaning is its purpose?.]
 			 	
But logic cannot address the natural history of the use of a word.
 			 	
To what extent can the use of one form of a word presuppose the existence of another form? For example, `dark' would presuppose `darker' and vice versa; or `white' would presuppose `whitish' and vice versa, etc.
\end{quote}
Ms-107,234 (1929--1930) (Ms-107: III, Philosophische Betrachtungen (\"ONB, Cod.Ser.n.22.020)) (Accessed 02 Jan 2025)
}
\end{quote}

Rather relevant to the point concerning the (Tractatus') view that our language was somehow an imperfect representation of the facts, here it comes another moment of self-correction which does show that there was not an abandonment of `language as a calculus' as such,\footnote{Here comes a remark from a manuscript of 1932-33 (Ms-114) saying:
\begin{quote}
Welche Rolle der Satz im Kalk\"ul spielt, das ist sein Sinn.
\end{quote}
which can be translated as:
\begin{quote}
What role the sentence plays in the calculus is its meaning.
\end{quote}
Ms-114,93r (Ms-114: X, Philosophische Grammatik (WL)) (Accessed 02 Jan 2025)
} but a recognition that `the mistake lies in not conceding to the calculus its real actual application, but in promising it for a nebulous distant ideal case':
\begin{quote}
Es ist von der gr\"o{\ss}ten Bedeutung, da\ss\ wir uns zu einem Kalk\"ul der Logik immer ein Beispiel seiner Anwendung denken, auf welches der Kalk\"ul wirklich eine Anwendung findet, \& da\ss\ wir nicht Beispiele, von denen wir || geben \& sagen, sie seien eigentlich nicht die idealen, diese aber h\"atten wir noch nicht. Das ist das Zeichen einer falschen Auffassung. (Russell \& ich haben, in verschiedener Weise an ihr laboriert. Vergleiche was ich in der ``Abhandlung" || ``Log. phil. Abh."  \"uber Elementars\"atze
\& Gegenst\"ande sage.) Kann ich den Kalk\"ul \"uberhaupt verwenden, dann ist dies auch die ideale Verwendung, \& die Verwendung um die es sich handelt. Einerseits will man n\"amlich das Beispiel nicht als das eigentliche anerkennen, weil man in ihm eine Mannigfaltigkeit sieht, der der Kalk\"ul nicht Rechnung tr\"agt. Anderseits ist es doch das Urbild des Kalk\"uls \& er davon hergenommen, \& auf eine getr\"aumte Anwendung kann man nicht warten. Man mu\ss\ sich also eingestehen, welches das eigentliche Urbild || Vorbild des Kalk\"uls ist.

     Nicht aber, als habe man damit einen Fehler begangen, den Kalk\"ul von daher genommen zu haben; sondern der Fehler || . Der Fehler liegt darin, dem Kalk\"ul seine wirkliche || eigentliche Anwendung jetzt nicht zuzugestehen \& sie || , sondern sie f\"ur eine nebulose Ferne || einen idealen Fall zu versprechen.\footnote{which can be translated as:
 \begin{quote}
 It is of the greatest importance that we always think of an example of its application to a calculus of logic, to which the calculus really finds an application, \& that we do not give examples of which we say they are not actually the ideal ones, but which we do not yet have. This is the sign of a false conception. (Russell \& I have laboured at it in different ways. Compare what I wrote in the `Tractatus' || `Log. phil. Abh.' on elementary propositions
\& objects). If I can use the calculus at all, then this is also the ideal use, \& the use in question. On the one hand, one does not want to recognise the example as the real one, because one sees in it a multiplicity that the calculus does not take into account. On the other hand, it is the archetype of the calculus \& it is taken from it, \& one cannot wait for a dreamed application. One must therefore admit to oneself which is the actual archetype of the calculus.

     Not, however, as if one had made a mistake in having taken the calculus from it; but the mistake lies in having taken the calculus from it. The mistake lies in not now conceding to the calculus its real actual application, but in promising it for a nebulous distant ideal case.
 \end{quote}
 Ms-115,55 (1931?) (Ms-115: XI, Philosophische Bemerkungen (WL)) (Accessed 26 December 2024)
 }
\end{quote}

With the recognition that logic at first would seem to serve as the `archetype of order', the writings in the so-called later phase redirected the attention to ordinary language, although this did not mean abandoning the view of language as a calculus, a (language-)game in which learning and understanding would be like being capable of seeing the role, the use and the purpose of words, signs, sentences in the game of language. With his own words, the intention was `to know what makes language language':
\begin{quote}
Immer wieder war man || ich versucht zu sagen ``Es mu\ss\ doch ?"
Wir haben in der Logik eine Theorie. Und die mu\ss\ einfach \& nett || ordentlich (neat) sein. Denn ich will ja wissen, wodurch die Sprache Sprache ist. Da\ss, was wir alles ``Sprache" nennen, Unvollkommenheiten, Schlacken, an sich hat, glaube ich, aber ich will das kennen lernen, was so verunreinigt ist. Das, wodurch ich im Stande bin, etwas zu sagen. Denn das mu\ss\ doch etwas sehr Eigent\"umliches sein. Der Gedanke, ein seltsames Wesen. (Er kann denken, was nicht der Fall ist.) Aber auch hier beruht das R\"atselhafte auf einem Mi{\ss}verst\"andnis || Mi{\ss}verst\"andnissen || Mi{\ss}verst\"andnis.

Ich habe mich ja seinerzeit gestr\"aubt gegen die Idee der nicht vollkommenen Ordnung in der Logik. ``Jeder Satz hat einen bestimmten Sinn"; ``In der Logik kann es nicht Unklarheit geben, denn sonst g\"abe es \"uberhaupt nicht Klarheit (\& also auch nicht Unklarheit).", ``Ein logisch-unklarer Satz ist || w\"are einer, der keinen bestimmten Sinn hat, also keinen Sinn". ? Hier spukt immer die Idee des \"atherischen Sinnes || Satzsinnes, dessen was man meint, des geistigen Prozesses.

     Die Logik schien das Urbild der Ordnung. Ich wollte immer (gegen Ramsey) sagen: Die Logik kann doch nicht zur empirischen Wissenschaft
werden. Aber wie wir die Sprache || W\"orter gebrauchen, das ist freilich Empirie || Erfahrung.

``Es l\"a{\ss}t sich eben doch denken!", oder ``Die Erfahrung zeigt || lehrt, da\ss\ es sich eben doch denken l\"a{\ss}t": das war solange revolting || emp\"orend, als Denken der geistige Proze\ss\ war || man im Denken den geistigen Proze\ss\ sah im Gegensatz zum Sprechen, Schreiben, etc. Die Logik mu{\ss}te solche Behauptungen entweder - dogmatisch - bestreiten, oder sich von ihnen zur\"uckziehen \& sagen da\ss\ sie da nichts zu tun habe, wo solche Fragen gestellt werden || beantwortet w\"urden || , wo auf solche Fragen geantwortet w\"urde. Aber wenn sie sich so zur\"uckzog, wo war || blieb dann noch ihr Feld? (Sie schien in || zu nichts zusammenzuschrumpfen.)\footnote{which can be translated as:
\begin{quote}
Again and again I was tempted to say ?It must be ...?
We have a theory in logic. And it has to be simple \& nice || neat (neat). Because I want to know what makes language language. I believe that everything we call `language' has imperfections, dross, in itself, but I want to get to know what is so polluted. That through which I am able to say something. For that must be something very peculiar. Thought, a strange creature. (It can think, which is not the case.) But here, too, the enigma is based on a misunderstanding || misunderstandings || misunderstanding.

At the time, I resisted the idea of imperfect order in logic. `Every sentence has a definite sense'; ``In logic there cannot be ambiguity, for otherwise there would be no clarity at all (\& therefore no ambiguity).", ``A logical-ambiguous sentence is || would be one that has no definite sense, therefore no sense". - Here the idea of the ethereal sense of the sentence, of what one means, of the mental process, always haunts.

     Logic seemed to be the archetype of order. I always wanted to say (against Ramsey): logic cannot become empirical science.
become an empirical science. But how we use language || words is, of course, empiricism || experience.

`It can be thought after all!'  or ``Experience shows || teaches that it can be thought after all": this was revolting || outrageous as long as thinking was the mental process || thinking was seen as the mental process in contrast to speaking, writing, etc.. Logic had to either - dogmatically - deny such assertions, or withdraw from them \& say that it had nothing to do where such questions were asked || answered || where such questions were answered. But if it withdrew like that, where was her field left? (It seemed to shrink to nothing in ||.)
\end{quote}
Ms-152,93 (1936) (Wittgenstein Nachlass Ms-152 [so-called C8] (WL) (Accessed 26 Dec 2024)
}
\end{quote}

In a manuscript of 1930 (Ms-183) there appears a remark to the effect of raising concerns on the attitude of `surrendering' to symbolism:
\begin{quote}
Ich glaube es geh\"ort heute Heroismus dazu die Dinge nicht als Symbole im Krausschen Sinn zu sehen. Das hei{\ss}t sich freizumachen von einer Symbolik, die zur Routine werden kann. Das hei{\ss}t freilich nicht versuchen sie wieder flach zu sehen sondern die Wolken des, sozusagen, billigen Symbolismus in einer h\"oheren Sph\"are wieder zu verdampfen (so da\ss\ die Luft wieder durchsichtig wird).

     Es ist schwer sich diesem Symbolismus heute nicht hinzugeben.
 	 			 	
     Mein Buch die log. phil. Abhandlung enth\"alt neben gutem \& echtem auch Kitsch d.h. Stellen mit denen ich L\"ucken ausgef\"ullt habe und sozusagen
in meinem eigenen Stil. Wie viele von dem Buch solche Stellen sind wei\ss\ ich nicht \& es ist schwer es jetzt gerecht zu sch\"atzen.\footnote{which can be translated as:
\begin{quote}
 I think it takes heroism today not to see things as symbols in the Krausian sense. That means freeing ourselves from a symbolism that can become routine. Of course, that doesn't mean trying to see them flat again, but to vaporise the clouds of, so to speak, cheap symbolism in a higher sphere (so that the air becomes transparent again).
 
     It is difficult not to give in to this symbolism today.
	 			 	
     My book the log. phil. Abhandlung contains not only good \& genuine but also kitsch, i.e. passages with which I have filled in gaps and, so to speak
in my own style. I don't know how many of the book are such passages \& it's hard to estimate it fairly now.
\end{quote}
Ms-183,30 (1930) (Ms-183 [so-called Denkbewegungen] (\"ONB, Cod.Ser.n.37.939)) (Accessed 02 Jan 2025)
}
\end{quote}
Indeed, in the \emph{Tractatus} one finds a `defense' of the use of (symbolic) logic in order to avoid `the most fundamental confusions (of which the whole of philosophy is full)', and that `we must employ a symbolism which excludes them':
\begin{quote}
3.323\\
In der Umgangssprache kommt es ungemein h\"aufig vor, dass dasselbe Wort auf verschiedene Art und Weise bezeichnet -- also verschiedenen Symbolen angeh\"ort --, oder, dass zwei W\"orter, die auf verschiedene Art und Weise bezeichnen, \"au{\ss}erlich in der gleichen Weise im Satz angewandt werden.\footnote{which was translated as:
\begin{quote}
3.323\\
In the language of everyday life it very often happens that the same word signifies in two different ways?and therefore belongs to two different symbols?or that two words, which signify in different ways, are apparently applied in the same way in the proposition. (Ogden)

3.323\\
In everyday language it very frequently happens that the same word has different modes of signification?and so belongs to different symbols?or that two words that have different modes of signification are employed in propositions in what is superficially the same way. (Pears \& McGuinness)
\end{quote}
}

3.324\\
So entstehen leicht die fundamentalsten Verwechslungen (deren die ganze Philosophie voll ist).\footnote{which was translated as:
\begin{quote}
3.324\\
Thus there easily arise the most fundamental confusions (of which the whole of philosophy is full). (Ogden)

3.324\\
In this way the most fundamental confusions are easily produced (the whole of philosophy is full of them). (Pears \& McGuinness)
\end{quote}
}

3.325\\
Um diesen Irrt\"umern zu entgehen, m\"ussen wir eine Zeichensprache verwenden, welche sie ausschlie{\ss}t, indem sie nicht das gleiche Zeichen in verschiedenen Symbolen, und Zeichen, welche auf verschiedene Art bezeichnen, nicht \"au{\ss}erlich auf die gleiche Art verwendet. Eine Zeichensprache also, die der logischen Grammatik -- der logischen Syntax -- gehorcht.

(Die Begriffsschrift Frege's und Russell's ist eine solche Sprache, die allerdings noch nicht alle Fehler ausschlie{\ss}t.)\footnote{which was translated as:
\begin{quote}
3.325\\
In order to avoid these errors, we must employ a symbolism which excludes them, by not applying the same sign in different symbols and by not applying signs in the same way which signify in different ways. A symbolism, that is to say, which obeys the rules of \emph{logical} grammar--of logical syntax.

(The logical symbolism of Frege and Russell is such a language, which, however, does still not exclude all errors.) (Ogden)

3.325\\
In order to avoid such errors we must make use of a sign-language that excludes them by not using the same sign for different symbols and by not using in a superficially similar way signs that have different modes of signification: that is to say, a sign-language that is governed by \emph{logical} grammar--by logical syntax.

(The conceptual notation of Frege and Russell is such a language, though, it is true, it fails to exclude all mistakes.) (Pears \& McGuinness)
\end{quote}
}

\end{quote}

In another Ms of 1932, we find an explicit reference to the view of `language as a calculus' which was not abandoned in the second phase, as it is claimed by some interpreters of Wittgenstein. But then, again, there comes an association of meaning of a proposition with its role, its purpose in the calculus of language, as well as with the consequences one can draw from it:
\begin{quote}
Der Sinn eines || des Satzes ist nicht pneumatisch, sondern ist das, was auf die Frage nach der Erkl\"arung des Sinnes zur Antwort kommt. Und -- oder -- der eine Sinn unterscheidet sich vom andern wie die Erkl\"arung des einen von der Erkl\"arung des andern.
			 	
     Welche Rolle der Satz im Kalk\"ul spielt, das ist sein Sinn.
		 	
     Der Sinn steht also nicht hinter ihm (wie der psychische Vorgang der Vorstellungen etc.).
 	 			 	
Welche S\"atze aus ihm folgen \& aus welchen S\"atzen er folgt, das macht seinen Sinn aus. Daher auch die Frage nach seiner Verifikation eine Frage nach seinem Sinn ist.\footnote{which can be translated as:
\begin{quote}
The sense of a proposition is not pneumatic, but is what comes to the answer to the question of the explanation of the sense. And -or - the one sense differs from the other as the explanation of the one differs from the explanation of the other.
			 	
     What role the proposition plays in the calculus is its meaning.
		 	
     The sense is therefore not behind it (like the mental process of ideas etc.).
 	 			 	
Which propositions follow from it \& from which propositions it follows, that constitutes its sense. Therefore the question of its verification is also a question of its meaning.
\end{quote}
Ms-113,42r (1932) (Ms-113: IX, Philosophische Grammatik (\"ONB, Cod.Ser.n.22.022)) (Accessed 09 Jan 2025)
}
\end{quote}

A few pages down in the same manuscript one encounters remarks in the spirit of self-criticism to the notion of correspondence between language and reality with a reference to the \emph{Tractatus}:
\begin{quote}
Nun kann man ruhig annehmen: ,ich meinte, Du solltest Licht machen' hei{\ss}t, da\ss\ mir dabei ein Phantasiebild von Dir in dieser T\"atigkeit vorgeschwebt hat, \& ebensogut, || : jener Satz || der Satz hei{\ss}t, da\ss\ mir dabei die Worte des vollst\"andigen Satzes in der Phantasie gegenw\"artig waren, oder da\ss\ eins von diesen beiden der Fall war, ? nur mu\ss\ ich wissen da\ss\ ich damit eine Festsetzung \"uber die Worte ,,ich meinte" getroffen habe \& eine engere als die ist welche dem tats\"achlichen allgemeinen Gebrauch des Ausdrucks entspricht.
			 	
Wenn das Meinen f\"ur uns irgendeine Bedeutung, Wichtigkeit, haben soll so mu\ss\ dem System der S\"atze ein System der Meinungen zugeordnet sein, was immer f\"ur Vorg\"ange die Meinungen sein sollen.

Inwiefern stimmt nun das Wort ,,Licht" im obigen Symbolismus oder Zeichenspiel mit einer Wirklichkeit \"uberein, ? oder nicht \"uberein?

      Wie gebrauchen wir \"uberhaupt das Wort \"ubereinstimmen? ? Wir sagen ,,die beiden Uhren stimmen \"uberein", wenn sie die gleiche Zeit zeigen, ,,die beiden Ma{\ss}st\"abe stimmen \"uberein", wenn gewisse Teilstriche zusammenfallen, ,,die beiden Farben stimmen \"uberein" wenn etwa ihre Zusammenstellung uns angenehm ist. Wir sagen ,,die beiden L\"angen stimmen \"uberein" wenn sie gleich sind, aber auch wenn sie in einem von uns gew\"unschten Verh\"altnis stehen. Und da\ss\ sie ,,\"ubereinstimmen" hei{\ss}t dann nichts andres, als da\ss\ sie in diesem Verh\"altnis - etwa 1 : 2 - stehen. So mu\ss\ also in jedem Fall erst festgesetzt werden, was unter ,,\"Ubereinstimmung" zu verstehen ist. - So ist es nun auch mit der

\"Ubereinstimmung einer L\"angenangabe mit einer L\"ange. Wenn ich sage: ,,dieser Stab ist 2 m lang", so kann ich z.B. erkl\"aren || eine Erkl\"arung geben, wie man nach diesem Satz mit einem Ma{ss}stab die L\"ange des Stabes kontrolliert, wie man etwa nach diesem Satz einen Me{\ss}streifen f\"ur den Stab erzeugt. Und ich sage nun der Satz stimmt mit der Wirklichkeit \"uberein, wenn der auf diese Weise konstruierte Me{\ss}streifen mit dem Stab \"ubereinstimmt. Diese Konstruktion eines Me{\ss}streifens illustriert \"ubrigens was ich in meinem Buch || der Abhandlung damit meinte da\ss\ der Satz bis an die Wirklichkeit herankommt. - Man k\"onnte sich das auch so klar machen: Wenn ich die Wirklichkeit daraufhin pr\"ufen will ob sie mit einem Satz \"ubereinstimmt, so kann ich das auch so machen da\ss\ ich sie nun beschreibe \& sehe ob der gleiche Satz herauskommt. Oder: ich kann die Wirklichkeit nach grammatischen Regeln in die Sprache des Satzes \"ubersetzen \& nun im Land der Sprache den Vergleich durchf\"uhren.

Als ich nun dem Andern erkl\"arte: ,,Licht" (indem ich Licht machte), ,,Dunkel || Finster" (indem ich ausl\"oschte) hatte ich auch sagen k\"onnen \& mit genau derselben Bedeutung ,,das ist || hei{\ss}t ,,Licht?" (wobei ich Licht mache) \& ,,das ist || hei{\ss}t ,,Finster" etc., und auch ebensogut: ,,das stimmt mit `Licht' \"uberein", ,,das stimmt mit ,Finster' \"uberein".
 
Es kommt eben wieder auf die Grammatik des Wortes ,,\"Ubereinstimmung" an, auf seinen Gebrauch. Und hier liegt die Verwechslung mit ,\"Ahnlichkeit' nahe in dem Sinn in welchem || dem zwei Personen einander \"ahnlich sind wenn ich sie leicht miteinander verwechseln kann.

     Ich kann auch wirklich nach der Aussage \"uber die Gestalt eines K\"orpers eine Hohlform konstruieren in die nun der K\"orper pa{\ss}t oder nicht pa{\ss}t, jenachdem die Beschreibung richtig oder falsch war; || , \& die konstruierte Hohlform geh\"ort dann in dieser Auffassung noch zur Sprache (die bis an die Wirklichkeit herankommt).\footnote{which can be translated as:
 \begin{quote}
 Now one may safely assume: `I meant that you should make light' means that I had an imaginary picture of you in this activity, \& just as well, || : that sentence || the sentence means that the words of the complete sentence were present to me in imagination, or that one of these two was the case, - only I must know that I have thereby made a determination about the words ?I meant? \& a narrower one than that which corresponds to the actual general use of the expression.
 
 If meaning is to have any significance, importance, for us, then a system of opinions must be assigned to the system of propositions, whatever the opinions are supposed to be.
	 			 	
     To what extent does the word ?light? in the above symbolism or play of signs correspond or not correspond to a reality?
     
      How do we even use the word coincide? - We say ?the two clocks agree? when they show the same time, ?the two scales agree? when certain graduation marks coincide, ?the two colours agree? when their combination is pleasant to us. We say ?the two lengths correspond? when they are the same, but also when they are in a relationship we desire. And that they ?match? then means nothing other than that they are in this ratio - about 1 : 2. So in each case we must first determine what is meant by ?correspondence?.

This is now also the case with the
Correspondence of a length specification with a length. If I say: ?This rod is 2 metres long?, I can explain, for example, how to check the length of the rod with a scale according to this theorem, how to create a measuring strip for the rod according to this theorem. And I now say that the theorem corresponds to reality if the measuring strip constructed in this way corresponds to the rod. By the way, this construction of a measuring strip illustrates what I meant in my book || the Tractatus by saying that the theorem comes close to reality. - It could also be explained like this: If I want to check reality to see whether it corresponds to a proposition, I can also do this by describing it and seeing whether the same proposition emerges. Or: I can translate the reality according to grammatical rules into the language of the sentence \& now carry out the comparison in the country of the language.

     When I now explained to the other: ?Light? (by making light), ?Dark || Darkness? (by extinguishing) I had also been able to say \& with exactly the same meaning
?this is || means `light' " (by making light) \& ``this is || means ``darkness" etc., and also just as well: ``this agrees with `light' ", ``this agrees with `darkness' ".
 \end{quote}
 Ms-113,46v (1932) (Ms-113: IX, Philosophische Grammatik (\"ONB, Cod.Ser.n.22.022)) (Accessed 02 Jan 2025)
 }
\end{quote}

\section{Late manuscripts with focus on meaning as \emph{Gebrauch}, \emph{Zweck}, \emph{Anwendung}, \emph{Verwendung}}
There are certainly a lot of repetitions in the \emph{Nachlass}, but in some cases the same remark may come in a different context, and with slight (but relevant) changes to the wordings of the sentences. In this section we are  interested in scrutinising the available material from the very late period, namely 1948--51, in order to get to the point where Wittgenstein goes `back to his roots' in the sense of making the connection between meaning, use, application, purpose of words, terms and sentences. To this end, it seems appropriate to take a look at the \emph{Nachlass} group entitled \emph{Collection 'Last Writings corpus' (1948-51)} by WAB Editors. The list of items of which this is part comes with a warning to the effect that: `[NB. Processing of the following may require a fair amount of time]'. (Accessed 05 Jan 2025). The Collection consists of Ms 137, 138, 169, 170, 171, 172, 173, 174, 175, 176, 177. 

For a start, it seems appropriate to bring in a series of remarks from Ms-173 which involve the words \emph{Bedeutung}, \emph{Anwendung}, \emph{Verwendung}, \emph{Gebrauch}, \emph{Konsequenzen}, \emph{Folgen}, in a rather intertwined fashion:
\begin{quote}
Denn auch wenn er selbst ohne zu l\"ugen sagt, ``er sei etwas \"argerlich gewesen, so hei{\ss}t das nicht, da\ss\ er damals in sich jenes von uns `\"argerlich' genanntes Gesicht gesehen habe. Wir haben wieder nur eine Wortreaktion von ihm, \& es ist noch gar nicht klar, wie viel die bedeutet. Das Bild ist klar, || ; aber nicht seine Anwendung.
	 			 	
     Denn auch, wenn ich selbst sage ``Ich habe mich etwas \"uber ihn ge\"argert"  -- wie wei\ss\ ich die Anwendung dieser Worte so genau? Ist sie denn so klar? Nun, sie sind eben eine \"Au{\ss}erung.
 	 			 	
     Aber wei\ss\ ich etwa nicht genau, was ich mit jener \"Au{\ss}erung meine? ``ich wei\ss\ doch genau, welchen Zustand in mir ich so nenne." Das hei{\ss}t nichts. Ich wei\ss, wie man das Wort anwendet \& manchmal mache ich die \"Au{\ss}erung ohne Z\"ogern || Bedenken \& manchmal z\"ogernd \& sage etwa, es war war nicht `geradezu \"argerlich' || ich h\"atte mich nicht `geradezu ge\"argert', oder dergleichen. Aber es ist nicht diese Unbestimmtheit, von der ich sprach. Auch dort, wo ich ohne Bedenken || unbedenklich sage, ich h\"atte mich ge\"argert, ist darum nicht ausgemacht wie sicher die Konsequenzen || weiteren Folgen aus diesem Signal sind.
	 			 	
     Als ich sagte es sei eine Unbestimmtheit in der Anwendung,
meinte ich nicht, ich wisse nicht recht, wann ich die \"Au{\ss}erung machen solle (wie es etwa w\"are, wenn ich nicht recht || gut Deutsch verst\"unde).
	 			 	
     Man mu\ss\ || darf eben nicht vergessen, welche Verbindungen gemacht werden, wenn wir lernen Ausdr\"ucke wie ``Ich \"argere mich" zu gebrauchen.
	 			 	
     Und denke nicht an ein Erraten der richtigen Bedeutung durch das Kind, denn, ob es sie richtig erraten hat, mu\ss\ sich doch wieder in seiner Verwendung der Worte zeigen.\footnote{which can be translated as:
\begin{quote}
For even if he himself says without lying that ?he was a little angry, that does not mean that he saw in himself the face we call ``angry?" Again we have only a verbal reaction from him, \& it is not yet clear how much it means. The image is clear, || ; but not its application.
 			 	
     For even if I myself say ?I was a little annoyed with him?, - how do I know the application of these words so exactly? Is it so clear? Well, they are just an utterance.

But don't I know exactly what I mean by that statement? `I know exactly what state of mind I'm talking about.?'That doesn't mean anything. I know how to use the word \& sometimes I make the utterance without hesitation || misgivings \& sometimes hesitantly \& say, for instance, it wasn't `downright angry?'|| I wouldn't have been `downright annoyed', or the like. But it's not that vagueness I was talking about. Even where I say without hesitation that I would have been annoyed, it is not clear how certain the consequences of this signal are.
 	 			 	
      When I said it was an indeterminacy in the application, I did not mean that I did not really know when I should make the statement (as it would be if I did not understand German very well).
      
We must not forget what connections are made when we learn to use expressions like `I'm annoyed'.
 	 			 	
     And don't think about the child guessing the correct meaning, because whether he has guessed it correctly must be shown again in his use of the words.
\end{quote}
Ms-173,46v (between 24 March and 12 April 1950) (Ms-173 [so-called Notebook no.5; source for Ms-176] (WL)) (Accessed 05 Jan 2025)
}
     
\end{quote}
Moving back from 1950, it may help to bring in a remark from \emph{Tractatus} connecting meaning and consequences which had appeared in a slightly different wording in the \emph{WW1 Notebooks 1914--16} dated 24.11.14:
\begin{quote}
5.1311\\
(Dass man aus $(x).fx$ auf $fa$ schlie{\ss}en kann, das zeigt, dass die Allgemeinheit auch im Symbol ,,$(x).fx$" vorhanden ist.)\footnote{which was translated as:
\begin{quote}
5.1311\\
(The fact that we can infer $fa$ from $(x).fx$ shows that generality is present also in the symbol ``$(x).fx$".) (Ogden)

5.1311\\
(The possibility of inference from $(x).fx$ to $fa$ shows that the symbol $(x).fx$ itself has generality in it.) (Pears \& McGuinness)
\end{quote}
}
\end{quote}
Fast forward to 1949, there comes an explicit association between meaning, use, application, consequences (\emph{Folgen}), purpose:
\begin{quote}
2.1.49.\\
     Was teilt uns der mit, der den gegenw\"artigen Aspekt meldet? sagt ``Jetzt seh ich es als ..." Das hei{\ss}t: Welche Folgen hat diese Mitteilung, welche Art von Verwendung kann sie haben? Sie k\"onnte verschiedenerlei Folgen haben.
     
     Wer z.B. den H.-E.-Kopf als H. sieht, wird nicht den Ausdruck des E.-Gesichts beschreiben k\"onnen.
     
     Raumvorstellung in der Darstellenden Geometrie. Wer das W\"urfelschema eben || jetzt flach sieht, wird verschiedene zeichnerische Operationen mit ihm nicht vornehmen k\"onnen. [Stimmt nicht ganz.]
	 			 	
     Verbindung mit dem Spiel ``Das k\"onnte ein ... sein".
 	 			 	
     Was teilst Du mir mit den || diesen || den Worten ... || mit dem Wort ... || durch die Worte mit? Was kann ich mit dieser \"Au{\ss}erung anfangen? Welche Folgen hat sie?\footnote{which can be translated as:
\begin{quote}
What does the person who reports the present aspect tell us? says `Now I see it as ...' This means: What are the consequences of this message, what kind of use can it have? It could have various consequences.
     For example, someone who sees the H. E. head as H. will not be able to describe the expression of the E. face.
     Spatial conception in descriptive geometry. If you see the cube scheme as flat, you will not be able to carry out various drawing operations with it. [Not quite true.]
	 			 	
     Connection with the game `This could be a ...?.
	 			 	
     What are you communicating to me with || these || the words ... || with the word ... || through the words? What can I do with this utterance? What consequences does it have?
\end{quote}
Ms-137,137b (1949) (Ms-137: R (WL)) (Accessed 08 Jan 2025)
}  

\end{quote}
Back to the early phase, there comes a series of remarks in the \emph{Tractatus} where meaning (\emph{Bedeutung}), sense (\emph{Sinn}), appear associated with \emph{Gebrauch}, \emph{Verwendung}, \emph{Folgen}:
\begin{quote}

3.326\\
Um das Symbol am Zeichen zu erkennen, muss man auf den sinnvollen Gebrauch achten.\footnote{which was translated as:
\begin{quote}
3.326\\
In order to recognize the symbol in the sign we must consider the significant use. (Ogden)

3.326\\
In order to recognize a symbol by its sign we must observe how it is used with a sense. (Pears \& McGuinness)
\end{quote}
}

3.327\\
Das Zeichen bestimmt erst mit seiner logisch-syntaktischen Verwendung zusammen eine logische Form.\footnote{which was translated as:
\begin{quote}
3.327\\
The sign determines a logical form only together with its logical syntactic application. (Ogden)

3.327\\
A sign does not determine a logical form unless it is taken together with its logico-syntactical employment. (Pears \& McGuinness)
\end{quote}
}

3.328\\
Wird ein Zeichen nicht gebraucht, so ist es bedeutungslos. Das ist der Sinn der Devise Occams.\footnote{which was translated as:
\begin{quote}
3.328\\
If a sign is \emph{not necessary} then it is meaningless. That is the meaning of Occam's razor. (Ogden)

3.328\\
If a sign is \emph{useless}, it is meaningless. That is the point of Occam's maxim. (Pears \& McGuinness)
\end{quote}
}
\end{quote}

Moving forward once more, it seems appropriate at this point to bring in a series of remarks from a very late manuscript where the interplay between the language of (symbolic) and ordinary language can shed light `on the borderline of logic and empiricism' and the importance of looking for the use/application (\emph{Verwendung}) when searching for the meaning, in this particular case, of words for colours:
\begin{quote}
Warum kann man sich wei{\ss}-durchsichtiges || durchsichtig-wei{\ss}es Glas nicht vorstellen, -- auch wenn es in Wirklichkeit keins gibt? Wo geht die Analogie mit dem durchsichtigen gef\"arbten || farbigen durchsichtigen schief?
	 			 	
     S\"atze werden oft an der Grenze von Logik \& Empirie gebraucht, (so) da\ss\ ihr Sinn \"uber die Grenzen || Grenze hin \& her wechselt, \& sie bald als Ausdruck einer Norm, bald als Ausdruck einer Erfahrung gelten.
     
     (Denn es ist ja nicht eine psychische Begleiterscheinung -- so stellt man sich den `Gedanken' vor --, sondern die Verwendung, die den logischen vom Erfahrungssatz unterscheidet.)\footnote{which can be translated as:
\begin{quote}
Why is it impossible to imagine white-transparent || transparent-white glass, - even if there is none in reality? Where does the analogy with the transparent coloured || coloured transparent go wrong?

Sentences are often used at the boundary of logic \& empiricism, (so) that their meaning changes back \& forth across the boundary, \& they are sometimes regarded as the expression of a norm, sometimes as the expression of an experience.

     (For it is not a mental concomitant - that is how one imagines the `thought' - but the use that distinguishes the logical from the empirical proposition).
\end{quote}
Ms-176,9r (21 March -- 24 April 1951) (Ms-176 [so-called Notebook no.6.1] (WL)) (Accessed 16 Jan 2025)
}
\end{quote}
Not much further down the Ms, there comes again a connection between the meaning of colour words to their use (\emph{Verwendung}):
\begin{quote}
Auf die Frage ``Was bedeuten die W\"orter `rot', `blau', `schwarz', `wei{\ss}' ", k\"onnen wir wohl || freilich gleich auf Dinge zeigen, die so gef\"arbt sind, -- aber weiter geht unsre F\"ahigkeit die Bedeutungen dieser Worte zu erkl\"aren nicht! Im \"ubrigen machen wir uns von ihrer Verwendung keine, oder eine ganz rohe, zum Teil falsche, Vorstellung.\footnote{which can be translated as:
\begin{quote}
In answer to the question ?What do the words ?red?, ?blue?, ?black?, ?white? mean??, we can certainly point to things that are coloured in this way - but our ability to explain the meanings of these words goes no further! For the rest, we have no, or a very crude, partly false, idea of their use.
\end{quote}
Ms-176,17r (21 March -- 24 April 1951) (Ms-176 [so-called Notebook no.6.1] (WL)) (Accessed 16 Jan 2025)
}
\end{quote}

Once again, a late manuscript such as Ms-175 (1950) brings up association between meaning, sense and application:
\begin{quote}
15.3.51\\
 ``I know that that's a tree." Warum kommt mir vor, ich verst\"unde den Satz nicht? obwohl || wo er || es doch ein h\"ochst || so einfacher Satz, von der aller gew\"ohnlichsten Art ist?
 
     Es ist als k\"onnte ich meinen Geist || ihn nicht auf irgend eine Bedeutung einstellen. Weil ich n\"amlich die Einstellung in einer Region || dem Bereiche suche, wo sie nicht || nicht in einer Region || dem Bereiche suche, wo sie ist. Sowie ich aus der philosophischen an die || eine allt\"agliche Anwendung dieses || des Satzes denke, wird sein Sinn klar \& gew\"ohnlich.
     
So wie die Worte ``Ich bin hier" nur in gewissen Zusammenh\"angen ? || Situationen ? Sinn haben, nicht aber, wenn ich sie Einem sage, der mir gegen\"uber sitzt \& mich klar sieht, -- \& zwar nicht darum, weil sie dann eine Selbstverst\"andlichkeit || \"uberfl\"ussig sind, sondern, weil ihr Sinn durch
die Situation nicht bestimmt ist, aber so eine Bestimmung braucht.\footnote{which can be translated as:
\begin{quote}
15.3.51\\
`I know that that's a tree.' Why does it seem to me that I do not understand the sentence? although it is a most simple sentence, of the most ordinary kind?

     It is as if I could not adjust my mind to any meaning. For I seek the attitude in a region where it is not, not in a region where it is.
     
As soon as I think from the philosophical to the everyday application of this sentence, its meaning becomes clear \& ordinary.
	 			 	
     Just as the words `I am here' only have meaning in certain contexts ? || situations ?, but not when I say them to someone who is sitting opposite me \& sees me clearly, - \& not because they are then self-evident || superfluous, but because their meaning is not determined by the situation, but thus needs a determination.
\end{quote}
Ms-175,49v (23 September 1950 -- 21 March 1951) (Ms-175 [so-called Notebook no.1] (WL)) (Accessed 16 Jan 2025)
}
\end{quote}

Back and forth from the early writings to the \emph{Investigations} and beyond helps us to see that a common thread in his view of language involves the connection between meaning, use, purpose, consequences. Even in a letter to Russell written in 1912 reveals such a connection in a rather explicit way:
\begin{quote}
Will you think that I have gone mad if I make the following
suggestion?: The sign $(x).\varphi x$ is not a complete
symbol but has meaning only in an inference of the kind: from $\vdash
\varphi x{\supset _ x}\psi x.\varphi (a)$
follows $\psi (a)$. Or more generally: from $\vdash (x).
\varphi x.{\varepsilon _ 0}(a)$ follows $\varphi (a)$. I
am---of course---most uncertain about the matter but something of the sort might
really be true.

(Letter {\it R.3\/} to Russell, dated 1.7.12, p.\ 12 of 
{\it Letters to Russell, Keynes and Moore\/})
\end{quote}

\section{Formal counterpart to purpose, use, consequences}
To the extent that one is willing to find a counterpart to meaning as use, purpose, consequences in formal logic, one needs to look at where does the `explanation of consequences' come in. Furthermore, since the language of formal logic does not get used in a direct `interaction' between an `Utterer' and an `Interpreter' (to borrow Peirce's terms), or, an `Attacker' and a `Defender' (Lorenzen), or even `Nature' and `Myself' (Hintikka), or `Spoiler' and `Duplicator' (Ehrenfecht-Fra\"{\i}ss\'e), in order to pinpoint where `purpose' comes into being, one needs to find where is this dichotomy to be found in the language of predicate logic and in the language of formalising calculable function (i.e., Church's lambda calculus). Coming to the rescue is, on the one hand, Church's notions of `normal form', abstraction-application, and term-reduction in $\lambda$-calculus, and, on the other hand, Gentzen's `introduction vs elimination' rules and proof reduction in Natural Deduction proof systems.\footnote{Just as a historical reminder, in both cases (predicate logic and lambda calculus) the actual languages come from the father of formal logic, namely, Gottlob Frege. Indeed, the language of predicate logic comes from his \emph{Begriffsschrift} and the language of the rules of calculating functions from his \emph{Grundgesetze} (with his `course-of-values' notation).}

In the context of the analysis of logical constants by rules of inference in Natural Deduction in the tradition of the verification principle, one finds an attempt to approach the formal counterpart of `purpose' in P.\ Martin-L\"of's discussion of `Correctness of assertion and validity of inference'\footnote{``This is a slightly edited transcript of a lecture given by Per Martin-L\"of on 26 October 2022 at the Rolf Schock Symposium in Stockholm."}:
\begin{quote}
As you can see from its title, my talk falls into two parts, and the order between
these two parts is essential: one has to begin with correctness of assertion and then
deal with validity of inference only after that. In addition to these two parts, I will
add a third part in the middle about the relation between correct assertion and
knowledge.
Let us begin, then, with the very notion of assertion, usually written like this,
$$\vdash C$$
Throughout, I use $C$ for the content of an assertion, and I use the Frege sign $\vdash$, or
the assertion sign in Russell?s terminology, for the assertoric force. An assertion in
general is just what you obtain by taking a content $C$ and prefixing the assertion
sign to it--a purely formal addition that you make.

(...)

The other part here, namely the content $C$, I have to say something about,
because it is given a semantic definition, simply as something to do. (...) As a term for something to do I will generally use
task, Aufgabe in German, which was introduced, albeit for propositions (Aussagen)
rather than assertoric contents, by Kolmogorov in 1932.

As far as the assertion sign is concerned, this explanation was a purely formal one.
From the contentual point of view, we also have to explain what is the purpose of
uttering an assertion. I take it that to explain the meaning of a complete sentence
is the same as explaining the purpose of the act of uttering it, that is, saying
something by means of it. So the question is, What is the purpose of an utterance
of the assertion $\vdash C$? To answer this question we have to introduce, not only the
speaker, who produces the assertion, but also the hearer, who receives the assertion
from the speaker, as indicated in the figure,
$$\mbox{speaker}\quad\longrightarrow\quad \vdash\ C \quad \longrightarrow\quad \mbox{hearer}$$
This is necessary, because we cannot explain the meaning, which is to say the
purpose, of an assertion speaking about the speaker alone: the meaning has to do
with the interaction between the speaker and the hearer. (Martin-L\"or 2022, p.1)\footnote{It is appropriate to point out here that this represents a major change in the position held by the proponents of the tradition (Heyting, Kolmogorov, Gentzen, Dummett, Prawitz, Martin-L\"of) of the so-called `verificationism' on which it is based a so-called language-based version of intuitionism:
\begin{quote}
So recall the explanations of the meanings of the
logical constants, the connectives and the quantifiers, given by Brouwer, Heyting and
Kolmogorov: they all follow the common pattern that, whatever the logical constant may be, an explanation is given of what a proof of a proposition formed by means
of that logical constant looks like, that is, what is the form, and, more precisely,
canonical or direct form, of a proof of a proposition which has that specific logical
constant as its outermost sign. It is clear from this what ought to be the general
explanation of what a proposition is, namely, that a proposition is defined by stipulating
how its proofs, more precisely, canonical or direct proofs, are formed. And, if we take
the rules by means of which the canonical proofs are formed to be the introduction
rules, I mean, if we call those rules introduction rules as Gentzen did, then his suggestion
that the logical constants are de fi ned by their introduction rules is entirely correct, so
we may rightly say that a proposition is de fi ned by its introduction rules.
\end{quote}
And, in the sequel:
\begin{quote}
Now what I would like to point out is that this is an explanation which could just
as well be identified with the verifi cation principle, provided that it is suitably interpreted. Remember first of all what the verification principle says, namely, that the
meaning of a proposition is the method of its verification. The trouble with that
principle, considered as a formula, or as a slogan, is that it admits of several different interpretations, so that there arises the question: how is it to be interpreted?
Actually, there are at least three natural interpretations of it. On the first of these, the
means of verifying a proposition are simply identified with the introduction rules for
it, and there is then nothing objectionable about Wittgenstein's formula, provided
that we either, as I just did, replace method by means, which is already plural in
form, or else make a change in it from the singular to the plural number: the meaning of a proposition is the methods of its verification. Interpreted in this way, it simply
coincides with the intuitionistic explanation of what a proposition is, or, if you
prefer, the Gentzen version of it in terms of introduction rules. For instance, using
this manner of speaking, there are two methods of verifying a disjunctive proposition,
namely, the two rules of disjunction introduction, and absurdity is defined by stipulating that it admits of no method of verification. (Martin-L\"of 2013, p.6)
\end{quote}
By the way, as we have remarked in previous opportunities, it is unwarranted to refer to `verification principle' as represented by conditions of proof as `Wittgenstein's formula', since he manifested a self-correction a couple of times to a previous remark he had made associating the meaning of a proposition to its method of proof.
}
\end{quote}
As we have emphasised in earlier occasions since the end of 1980's, the explanation of (immediate) consequences of a proposition in the context of Natural Deduction rules is made explicit by the so-called reduction rules which spell the effect of the elimination rules on the result of the application of introduction rules. And this opens the way to see a parallel between the pragmatist/dialogical/game-theoretic accounts of the semantics of predicate logic: Utterer/ Defender/ Myself come in via the introduction rules whereas Interpreter/Attacker/Nature find their counterpart in the elimination rules. And the Inversion Principle, in the form of reduction rules, brings them together for logical harmony.

\subsection{Lambda calculus and the rule of application ($\beta$-reduction)}

As far as the semantics of Alonzo Church's $\lambda$-calculus, which he introduces in (1936) `to propose a definition of effective calculability', the actual meaning-giving rule is the rule of application (so-called $\beta$-reduction) which essentially describes, via a process of substitution, the role and the purpose of the $\lambda$-symbol binding a variable and working as an `abstraction' operator, cf.\ Rule II below:
\begin{quote}
I. To replace any part $\lambda x[M]$ of a formula by $\lambda y[{\mathrm S}^x_y M |]$, where $y$ is
a variable which does not occur in $M$.

II. To replace any part $\{\lambda x[M]\}(N)$ of a formula by ${\mathrm S}^x_N M | $, provided
that the bound variables in $M$ are distinct both from $x$ and from the free
variables in $N$. 

III. To replace any part ${\mathrm S}^x_N M | $ (not immediately following $\lambda$) of a formula by $\{\lambda x [M]\}  (N)$, provided that the bound variables in $M$ are distinct
both from $x$ and from the free variables in $N$. (Church 1936, p.347)
\end{quote}
Rules I and III were given by Church to, respectively, formalise the change in bound variable and to give the inverse of rule II in order to make the conversion relation (sequences of applications of rules I, II and III) an equivalence relation.
This way, a $\lambda$-term represents a function `as a rule', rather than the usual set-theoretic representation as a set of tuples whose last component is the result of applying the function to the previous components of the tuple.\footnote{
As a matter of historical as well as technical record, it is appropriate to bring here the work of Frege in the formalisation of `function as a rule' via a notation which finds direct parallels to Church's. In spite of a clear emphasis on formulating precisely the notion of the `extension' of a concept, Frege's {\em Grundgesetze\/} calculus of expressions of the form `$\acute{e}f(e)$', which goes back to his (1891) `Function and Concept' \emph{Werthverlauf}, brought up a number of operational aspects of what today is seen as the theory of `functions as rules' as defined by Church with his $\lambda$-calculus.
\begin{quote}
Wenn ich allgemein sage:

$\qquad$ es bedeute

$$,\acute{\varepsilon}\Phi(\varepsilon)\mbox{`}$$

$\qquad$ den Werthverlauf der Function $\Phi(\xi)$

so bedarf dies ebenso einer Erg\"angzung wie oben unsere Erkl\"arung von ,$\forall \alpha \Phi(\alpha)$` (\emph{Grundgesetze} I, \S 9, p. 15)
\end{quote}
which was translated as
\begin{quote}
If I say in general:

$\qquad$ let

$$\mbox{`}\acute{\varepsilon}\Phi(\varepsilon)\mbox{'}$$

$\qquad$ refer to the value-range of the function $\Phi(\xi)$,

then this too requires supplementation, just like our explanation of $\forall \alpha \Phi(\alpha)$
above.
\end{quote}
The device of variable-binding, and the idea of having terms representing
incomplete `objects' whenever they contain free variables, were both introduced
in a systematic way by Frege in his {\it Grundgesetze\/}. As early as 1893 Frege
developed in his {\it Grundgesetze I\/} what can be seen as the early origins
of the notions of {\it abstraction\/} and {\it application\/}, when showing
techniques for transforming functions (expressions with free variables) into
value-range terms (expressions with no free variables) by means of an
`introductory' operator of abstraction producing the {\it Werthverlauf\/}
expression, e.g., `$\acute{\varepsilon}f(\varepsilon)$', and the effect of its
corresponding `eliminatory' operator `$\scriptscriptstyle \cap$' on a
{\it value-range\/} expression. 

Expressing how important he considered the introduction of a variable-binding
device for the functional calculus (recall that the variable-binding device for
the logical calculus (i.e., the universal quantifier) had been introduced earlier in {\it Begriffsschrift\/}),
Frege says:
\begin{quote}
The introduction of a notation for value-ranges seems to
me to be one of the most consequential additions to my concept-script that I made
since my first publication on this matter.
\hfill{({\it Grundgesetze I\/}, \S 9, p. 15f.)}
\end{quote}
In order to make sure that a value-range term (or the term resulting of an application of this to another term) would not end up without a denotation, Frege introduces the operator `$\backslash$':
\begin{quote}
(...) in the case in which $\Phi(\xi)$ was a concept under which one and only one object fell, ``$\acute{\varepsilon}\Phi(\varepsilon)$'' would designate this object. Now of course this is not possible, because that equation could not be sustained in its general form, but we can serve our purpose by introducing the function
$$\backslash\xi$$
with the specification to distinguish two cases:\\
1. If, for the argument, there is an object $\Delta$ such that $\acute{\varepsilon}(\Delta=\varepsilon)$ is the argument, then the value of the function $\backslash\xi$ is to be $\Delta$ itself;\\
2. if, for the argument, there is no object $\Delta$ such that $\acute{\varepsilon}(\Delta=\varepsilon)$ is the argument, then the argument itself is to be the value of the function $\backslash\xi$. (\S11, p. 19)
\end{quote}
By doing this, Frege makes sure that ``$\backslash\acute{\varepsilon}\Phi(\varepsilon)$'' always has a reference. Now, it may be worth noticing the similarity between the two notations
$$\mbox{Frege's }\backslash\acute{\varepsilon}\Phi(\varepsilon)\qquad\mbox{and}\qquad \mbox{Church's }\lambda \varepsilon\Phi(\varepsilon)$$
which are both terms representing a function as a rule. Of course, in Church's case, there was no need to make sure the term always had a reference since the $\lambda$-calculus was meant to formalise {\em partial\/} functions.

In order to define the use, the purpose of value-range terms, Frege defines the operation of function-application as (Vol.\ I, \S 34, pages 52ff), (translated by M.\ Furth):
\begin{quote}
``(...) it is a matter only of designating the value of the function $\Phi(\xi)$
for the argument $\Delta$, i.e., $\Phi (\Delta)$, by means of ``$\Delta$'' and
``$\acute{\varepsilon}\Phi(\varepsilon)$''. I do so in this way:
$$\mbox{``}\Delta\ {\scriptscriptstyle \cap}\ \acute{\varepsilon}\Phi(\varepsilon)\mbox{''}$$
which is to mean the same as ``$\Phi(\Delta)$''.''
\end{quote}
\noindent Note the similarity to the rule of functional {\it application\/} (Rule II of Church 1936),
where `$\Delta$' is the argument, `$\acute{\varepsilon}\Phi(\varepsilon)$'
is the function, and `$\scriptscriptstyle \cap$' is the application operator.

Having settled the notation and the meaning of the symbol for the formalisation of `function as a rule', Frege goes on to show how to use this formal apparatus to define `numeral' (as opposed to `number') as a term in the calculus containing value-ranges, much in the same way Church shows how to formalise effective calculability in Arithmetic in his $\lambda$-calculus.
}
 Arising out of Church's work of formulating mathematically the notion of `computable function' in order to prove that Hilbert's Entscheidungsproblem was unsolvable, the rules of the calculus would not prevent a function from being applied to itself, so that self-application would be `legal' and this led to the first formalisation of the notion of `partial function' in the literature. But this meant that the `set-theoretic' meaning of this formal counterpart to partial function had to be elaborated in a more sophisticated way leading to substantial work on what is called denotational semantics, out of which Dana Scott's domain theory is the most well known. This way, the search for a mathematical `understanding' of a $\lambda$-term has taken two main directions: operational semantics and denotational semantics, the latter being an approach to the `meaning' of a $\lambda$-term via an interpretation into the universe of (ordered, topologically structured) sets. As to the former, the foundational work was done by G.\ Plotkin (1974) which helped to consolidate the theory of programming languages, even if the first use of the term `operational semantics' is credited to D.\ Scott in a seminal paper in 1970. In any case, it was Church/Frege who prepared the grounds for the study of `function as a rule' by stating the rule of reduction of a $\lambda$-term via making explicit its meaning as use, purpose, consequences, i.e., as a general rule which is to be applied to any input resulting in a new term where such input replaces the bound variable occurring in the original term. That is the most basic settlement of meaning for the $\lambda$-calculus, which comes first and foremost before any (most legitimate!) approach to associate a term to a function between natural numbers, be it via the application of a sequence of rules of reduction, be it via an interpretation into the realm of (ordered, topological) sets.

\subsection{Predicated logic and rules of (proof-)reduction}
Coming out of H.\ Curry's (1934) discovery of parallels between the axioms of (intuitionistic) implicational logic and the types of the combinators of Combinatory Logic, extended by W.\ Howard's (1969) to full (intuitionistic) predicate logic, there is a formalisation of Errett Bishop's style of constructive mathematics which was then called intuitionistic type theory by its main author, Per Martin-L\"of. The idea behind what is called `Curry-Howard interpretation of predicate logic leading to the slogan `propositions are types (of their proofs)' is to associate typing rules to proof rules in the style of Natural Deduction. And this association finds a direct counterpart to the so-called Brouwer-Heyting-Kolmogorov (BHK) interpretation of intuitionistic logic, where proofs, not truth-values, determine the meaning of propositions.
$$\begin{array}{rcll}
\bot & \mbox{as} &\varnothing&\mbox{(empty type)}\\
\alpha\land\beta& \mbox{as}& \alpha\times\beta&\mbox{(product)}\\ 
\alpha\lor\beta&\mbox{as}&\alpha+\beta&\mbox{(sum)}\\
\forall x\gamma&\mbox{as}&\Pi x\gamma &\mbox{(dependent product) }\\
\exists x\gamma&\mbox{as}&\Sigma x\gamma &\mbox{(dependent sum) }\
\end{array}$$
Here is where one can start to look at where meaning as use, purpose, consequences is supposed to show up. For a starter, it may be easier to harness on the previous discussion of $\lambda$-terms, given that, not just type theory is considered to be a formalisation of typed $\lambda$-calculus, but, according to the BHK-interpretation, a proof of an implication $A\to B$ is a function which transforms proofs of $A$ into proofs of $B$. In other words, formal proofs of $A\to B$ are $\lambda$-terms. And we have discussed above where the `meaning as use, purpose' of $\lambda$-terms become clear: in the $\beta$-reduction rules. By `lifting' this up to the type-rules for $A\to B$, we have:

\medskip

$\rightarrow$-$\beta$-{\it reduction\/}
$$\displaystyle{{\displaystyle{\ \atop {a:A}} \qquad \displaystyle{{\displaystyle{{[x:A]} \atop {b(x):B}}} \over
{\lambda x.b(x):A\rightarrow B}}\rightarrow\mbox{\it -intr}} \over
{{\tt APP}(\lambda x.b(x),a):B}}\rightarrow\mbox{\it -elim} \quad \beta\mbox{-reduces to}
\quad \displaystyle{{a:A} \atop {b(a/x):B}}$$
Associated rewriting:\\
${\tt APP}(\lambda x.b(x),a)=_\beta b(a/x)$

\medskip

Taking from the rules of \emph{reduction} between proofs in a system of Natural Deduction enriched with terms alongside formulas, such as in the so-called Curry--Howard interpretation (Howard (1980)) of which Martin-L\"of's type theory can be seen as an instance, and isolating the terms corresponding the derivations, one can draw the following parallels:
Looking at the conclusion of reduction inference rules, one can take the destructor as being the Attack (or `Nature', its counterpart in the terminology of Hintikka's Game-Theoretical Semantics, or `Interpreter' as in Peirce's formulation), and the constructor as being the Assertion (or Hintikka's `Myself', or Peirce's `Utterer'). This way the game-theoretic explanations of logical connectives find direct counterpart in the functional interpretation with the view that meaning is use, purpose, application, consequences:

\medskip

\noindent $\land$-$\beta$-{\it reduction\/}
$$\displaystyle{{\displaystyle{{a:A \qquad b:B} \over
{\langle a,b\rangle:A\land B}}\land\mbox{\it -intr}} \over
{{\tt FST}(\langle a,b\rangle):A}}\land\mbox{\it -elim} \qquad \qquad
\beta\mbox{-reduces to} \qquad \qquad a:A$$
$$\displaystyle{{\displaystyle{{a:A \qquad b:B} \over
{\langle a,b\rangle:A\land B}}\land\mbox{\it -intr}} \over
{{\tt SND}(\langle a,b\rangle):B}}\land\mbox{\it -elim} \qquad \qquad
\beta\mbox{-reduces to} \qquad \qquad b:B$$
Associated rewritings:\\
${\tt FST}(\langle a,b\rangle)=_\beta a$\\
${\tt SND}(\langle a,b\rangle)=_\beta b$\\
\ 
\\
\begin{tabbing}
\underline{Assertion}/{\it Introd.\/} $\quad$ \= \underline{Attack}/{\it Elim.\/} $\qquad\qquad\quad \beta\mbox{-reduces to} \qquad\qquad \qquad
$ \= \underline{Defense} \\
\\
Conjunction (`$\land$'): \\
${A}\land{B}$ \> $\quad \quad$ $L$? \> ${A}$ \\
${A}\land{B}$ \> $\quad \quad$ $R$? \> ${B}$ \\
\\
$\langle a,b\rangle:{A}\land{B}$ \> {\tt FST}$(\langle a,b\rangle)$ $\qquad \beta\mbox{-reduces to}$ \> $a:{A}$ \\
$\langle a,b\rangle:{A}\land{B}$ \> {\tt SND}$(\langle a,b\rangle)$ $\qquad \beta\mbox{-reduces to}$ \> $b:{B}$ \\
\end{tabbing}

\noindent $\lor$-$\beta$-{\it reduction\/}
$$\displaystyle{{\displaystyle{{a:A} \over {{\tt inl}(a):A\lor B}}\lor\mbox{\it -intr} \  \displaystyle{{[x:A]} \atop {f(x):C}} \  \displaystyle{{[y:B]} \atop {g(y):C}}} \over
{{\tt CASE}({\tt inl}(a),\upsilon x.f(x),\upsilon y.g(y)):C}}\lor\mbox{\it -elim} \ \ \beta\mbox{-reduces to} \  \displaystyle{{a:A} \atop {f(a/x):C}}$$
$$\displaystyle{{\displaystyle{{b:B} \over {{\tt inr}(b):A\lor B}}\lor\mbox{\small\it -intr} \  \displaystyle{{[x:A]} \atop {f(x):C}} \  \displaystyle{{[y:B]} \atop {g(y):C}}} \over
{{\tt CASE}({\tt inr}(b),\upsilon x.f(x),\upsilon y.g(y)):C}}\lor\mbox{\small\it -elim} \ \ \beta\mbox{-reduces to} \ \  \displaystyle{{b:B} \atop {g(b/y):C}}$$
Associated rewritings:\\
${\tt CASE}({\tt inl}(a),\upsilon x.f(x),\upsilon y.g(y))=_\beta f(a/x)$\\
${\tt CASE}({\tt inr}(b),\upsilon x.f(x),\upsilon y.g(y))=_\beta g(b/y)$\\
\ 
\\
\begin{tabbing}
\underline{Assertion}/{\it Introd.\/} $\quad$ \= \underline{Attack}/{\it Elim.\/} $\qquad\qquad\quad \beta\mbox{-reduces to} \qquad\qquad \qquad
$ \= \underline{Defense} \\
\\

Disjunction (`$\lor$'): \\
${A}\lor{B}$ \> $\quad \quad$ ? \> ${A}$ \\
${A}\lor{B}$ \> $\quad \quad$ ? \> ${B}$ \\
\\
{\tt inl}$(a):{A}\lor{B}$ \> {\tt CASE}({\tt inl}$(a),\upsilon x.f(x),\upsilon y.g(y))$ $\quad \beta\mbox{-reduces to}$ \> $a:A$, $f(a/x):{C}$ \\
{\tt inr}$(b):{A}\lor{B}$ \> {\tt CASE}({\tt inr}$(b),\upsilon x.f(x),\upsilon y.g(y))$ $\quad \beta\mbox{-reduces to}$ \> $b:B$, $g(b/y):{C}$ \\
\\
\end{tabbing}
Whilst in both {\tt FST}$(\langle a,b\rangle)$ and {\tt SND}$(\langle a,b\rangle)$ the destructors allow access to either of the conjuncts, in {\tt CASE}({\tt inl}$(a),\upsilon x.f(x),\upsilon y.g(y))$ and
 {\tt CASE}({\tt inr}$(b),\upsilon x.f(x),\upsilon y.g(y))$ the destructor is not given access to the either disjunct but must ask for whichever disjunct comes from the introduction.

\medskip

\noindent $\rightarrow$-$\beta$-{\it reduction\/}
$$\displaystyle{{\displaystyle{\ \atop {a:A}} \qquad \displaystyle{{\displaystyle{{[x:A]} \atop {b(x):B}}} \over
{\lambda x.b(x):A\rightarrow B}}\rightarrow\mbox{\it -intr}} \over
{{\tt APP}(\lambda x.b(x),a):B}}\rightarrow\mbox{\it -elim} \quad \beta\mbox{-reduces to}
\quad \displaystyle{{a:A} \atop {b(a/x):B}}$$
Associated rewriting:\\
${\tt APP}(\lambda x.b(x),a)=_\beta b(a/x)$\\
\ 
\\
\begin{tabbing}
\underline{Assertion}/{\it Introd.\/} $\quad$ \= \underline{Attack}/{\it Elim.\/} $\qquad\qquad\quad \beta\mbox{-reduces to} \qquad\qquad \qquad
$ \= \underline{Defense} \\
\\
Implication (`$\rightarrow$'): \\
${A}\rightarrow{B}$ \> $A$ $\quad$ ? \> {$B$} \\
\\
$a:A$\\
$\lambda x.b(x):A\to B$ \> {\tt APP}$(\lambda x.b(x),a)$ $\qquad
\beta\mbox{-reduces to}$ \> $b(a/x):{B}$ \\
\\
\end{tabbing}

\noindent $\forall$-$\beta$-{\it reduction\/}
$$\displaystyle{{\displaystyle{\ \atop {t:D}} \qquad \displaystyle{{\displaystyle{{[x:D]} \atop {g(x):P(x)}}} \over
{\Lambda x.g(x):\forall x^D.P(x)}}\forall\mbox{\it -intr}} \over
{{\tt EXTR}(\Lambda x.g(x),t):P(t)}}\forall\mbox{\it -elim} \quad \beta\mbox{-reduces to}
\quad \displaystyle{{a:D} \atop {g(t/x):P(t)}}$$
Associated rewriting:\\
${\tt EXTR}(\Lambda x.g(x),t)=_\beta g(t/x)$\\
\ 

Here it may be helpful to quote Wittgenstein on the meaning of `all' (and of `$\forall x P(x)$', for that matter):
\begin{quote}
Man k\"onnte sagen: Man lernt die Bedeutung von ``alle", indem man lernt, da\ss\ aus $(x).fx$ $fa$ folgt. -- D.h., die \"Ubungen die den Gebrauch dieses Wortes ein\"uben, lehren, laufen || gehen immer darauf hinaus, da\ss\ keine Ausnahme gemacht werden darf. || , die den Gebrauch dieses Wortes ein\"uben, -- seine Bedeutung lehren, || -- zielen immer dahin, da\ss\ eine Ausnahme nicht gemacht werden darf.

     ``Aus `alle', wenn es so gemeint ist mu\ss\ doch das folgen." --Wenn es wie gemeint ist? \"Uberlege es Dir, wie meinst Du es? Da schwebt Dir etwa noch ein Bild vor -- \& mehr hast Du nicht. -- Nein, es mu\ss\ nicht, -- aber es folgt: Wir vollziehen diesen \"Ubergang.
     
     Und wir sagen: Wenn es || das nicht folgt, dann waren es eben nicht alle! -- -- und das zeigt nur, wie wir mit Worten in so einer Situation reagieren. --
     
     Wir k\"onnten es auch so sagen: Es kommt uns vor, da\ss, wenn aus $(x). fx$ nicht mehr $fa$ folgen soll, au{\ss}er dem Gebrauch des Wortes ``alle" noch etwas anderes sich ge\"andert hat || haben mu\ss, etwas, was dem Worte unmittelbar || selbst anh\"angt.\footnote{which can be translated as:
\begin{quote}
You could say: You learn the meaning of `all' by learning that $fa$ follows from $(x).fx$. - In other words, the exercises that practise, teach, run || the use of this word always result in the fact that no exception may be made. || , which practise the use of this word, - teach its meaning, || - always aim at the fact that an exception must not be made.

     `From `all', if it is meant that way, this must follow.' - If it is meant how? Think about it, what do you mean? You have another image in mind - \& that's all you have. - No, it doesn't have to - but it follows: We make this transition.
     
     And we say: If it doesn't follow, then it wasn't all of them! - and that just shows how we react with words in such a situation. -
     
     We could also put it this way: It seems to us that if fa is no longer to follow from $(x).fx$, something else must have changed besides the use of the word `all', something that is directly attached to the word itself.
 \end{quote}
Ms-117,12 (1937) (Ms-117: XIII, Philosophische Bemerkungen (WL)) (Accessed 16 Jaan 2025)
 }
\end{quote}
In the same Ms, we find a remark in the same spirit:
\begin{quote}
Und man macht den Sinn von ``$(x).fx$" klar, indem man darauf dringt, da\ss\ aus ihm ``$fa$" folgt.\footnote{which can be translated as:
\begin{quote}
And we make the meaning of `$(x).fx$' clear by insisting that `$fa$' follows from it.
\end{quote}
Ms-117,97(1937) (Ms-117: XIII, Philosophische Bemerkungen (WL)) (Accessed 16 Jaan 2025)

\medskip

A similar remark from 24.11.14 which was quoted above already appears in the \emph{WW1 Notebooks 1914--16} (Ms-102,37r)
}
\end{quote}
As we can observe from the rules of the existential quantifier (below), one cannot explain the immediate consequences of a sentence such as $\exists x^D.P(x)$ by presenting a specific witness, since it does not state that any witness can be used to instantiate the predicate $P(-)$, and that is why the rule of $\exists$-\emph{elimination} has to bring in a new local assumption of the kind `suppose the witness is $t$', and then somehow get rid of the special status of such witness. And special conditions apply when handling the deduction from such a new (local) assumption to the final conclusion $C$. In the dialogical/game-theoretical approaches such issue is dealt with the use of the dichotomy Attacker/Defender, Nature/Myself, Interpreter/Utterer and the shift of the choice of the witness to the Attacker (or Nature, or Interpreter).

\medskip

\noindent $\exists$-$\beta$-{\it reduction\/}
$$\displaystyle{{\displaystyle{{s:D \quad g(s):P(s)} \over
{\varepsilon x.(g(x),s):\exists x^D.P(x)}}\exists\mbox{\it -intr} \quad
\displaystyle{{[t:D,h(t):P(t)]} \atop {d(h,t):C}}} \over
{{\tt INST}(\varepsilon x.(g(x),s),\sigma h.\sigma t.d(h,t)):C}}\exists\mbox{\it -elim} \  \beta\mbox{-reduces to}\ 
\displaystyle{{s:D,g(s):P(s)} \atop {d(g/h,s/t):C}}$$
Associated rewriting:\\
${\tt INST}(\varepsilon x.(g(x),s),\sigma h.\sigma t.d(h,t))=_\beta d(g/h,s/t)$\\
\

\begin{tabbing}
\underline{Assertion}/{\it Introd.\/} $\quad$ \= \underline{Attack}/{\it Elim.\/} $\qquad\qquad\quad \beta\mbox{-reduces to} \qquad\qquad \qquad
$ \= \underline{Defense} \\
\\
Universal Quantifier (`$\forall$'): \\
$\forall x^{D}.{P}(x)$ \> $a:{D}$ ? \> ${P}(a)$ \\
\\
$a:A$\\
$\Lambda x.g(x):\forall x^{D}.{P}(x)$ \\
\> {\tt EXTR}$(\Lambda x.g(x),t)$ $\qquad
\beta\mbox{-reduces to}$ \> $g(t/x):{P}(t)$ \\
\\
Existential Quantifier (`$\exists$'): \\
$\exists x^{D}.{P}(x)$ \> $\quad \quad$ ? \> $s:{D}, {P}(s)$ \\
\\
\> \> $s:D, g(s):P(s)$\\
$\varepsilon x.(g(x),s):\exists x^{D}.{P}(x)$ \\
\> {\tt INST}$(\varepsilon x.(g(x),s),\sigma h.\sigma t.d(h(t),t))$ $\ \ \beta\mbox{-reduces to}$ \> $d(g(s)/h(t),s/t):{C}$ \\
\\
\end{tabbing}
Whilst in {\tt EXTR}$(\Lambda x.g(x),t)$  the destructor has access to the witness (can use a generic element)
$a$, in {\tt INST}$(\varepsilon x.(g(x),s),\sigma h.\sigma t.d(h(t),t))$ the only option is to eliminate over the witness $s$ which is inside the term $\varepsilon x.(g (x), s)$ built with the constructor because it was chosen by the introduction.

\medskip

\noindent $Id$-$\beta$-{\it reduction\/}
$$\displaystyle{{\displaystyle{{a=_r b:A} \over {r(a,b):Id_A(a,b)}}Id\mbox{\it -intr\/} \qquad
\displaystyle{{[a=_t b:A]} \atop {d(t):C}}} \over
{{\tt REWR}(r(a,b),\sigma t.d(t)):C}}Id\mbox{\it -elim\/} \qquad
\beta\mbox{-reduces to} \qquad
\displaystyle{{a=_r b:A} \atop {d(r/t):C}}$$
Associated rewriting:\\
${\tt REWR}(r(a,b),\sigma t.d(t))=_\beta d(r/t)$\\
\
\begin{tabbing}
\underline{Assertion}/{\it Introd.\/} $\quad$ \= \underline{Attack}/{\it Elim.\/} $\qquad\qquad\quad \beta\mbox{-reduces to} \qquad\qquad \qquad
$ \= \underline{Defense} \\
\\
Propositional Equality (`$Id_A(a,b)$'): \\
$Id_A(a,b)$ \> $\quad \quad$ ? \> $a=_r b:A$ \\
\\
$r(a,b):Id_A(a,b)$ \> {\tt REWR}$(r(a,b),\sigma t.d(t))$ $\ \ \beta\mbox{-reduces to}$ \> $d(r/t):{C}$ \\

\end{tabbing}

In {\tt REWR}$(r(a,b),\sigma t.d(t))$ the destructor has no choice regarding the 
`reason' for $a$ being equal to $b$, since $r$ will have been chosen by the time of the assertion, i.e.\ the application of the introduction rule.

\medskip

The reference to these parallels is intended to exhibit a connection with the ``pragmatist"/``dialogical" approaches to meaning as revealed in Peirce's writings on the interaction between the \emph{Interpreter} and the \emph{Utterer} which seem to pertain to the ``game"/``dialogue" approaches to meaning (Lorenzen, Hintikka, Ehrenfeucht-Fra\"{\i}ss\'e).

In any case, it seems fair to say that it is due to the pioneers of the introduction of `dialogical/game-theoretic/pragmatic' aspects, in the names of Peirce, Lorenzen, Hintikka, Ehrenfeucht-Fra\'{\i}ss\'e, who prepared the grounds for the study of `semantics' of predicate logic by stating the rules of `attack-defense'/`Nature-Myself'/`Utterer-Interpreter'/elimination-introduction (reduction) for each logical connective and quantifier of a formula via making explicit its meaning as use, purpose, consequences, i.e., as a general rule which is to translate into precise `rules of interaction' when it comes to sentences and their meaning, use, purpose, consequences. That would be essentially the most basic settlement of meaning for predicate logic, which comes first and foremost before any (most legitimate!) approach to associate a sentence to a relation of satisfaction within a model in the realm of set theory..

\section{Conclusions}
By looking at Wittgenstein's \emph{Nachlass} at large, we have gathered substantial evidence to add support to the view that there is a common thread in his thinking towards the view that meaning is use, purpose, application, usefulness, in spite of twists and turns of a restless mind. Rather relevant to the point concerning the (Tractatus') view that our language was somehow an imperfect representation of the facts, one can find several moments of self-correction which reveal that there was not a straight rejection of the view of `language as a calculus' as such. With the recognition that logic at first would seem to serve as the `archetype of order', as it happened in the picture theory of the Tractatus, the writings in the so-called later phase redirected the attention to ordinary language, although this did not mean moving away from the conception of language as a calculus, a (language-)game in which learning and understanding would be like being capable of seeing the role, the use and the purpose of words, signs, sentences in the game of language. With his own words, the intention was `to know what makes language language'.

With the same intention as in previous papers, we have made an incursion into the origins of Church's $\lambda$-calculus and Frege's Predicate Logic in order to emphasise the fact that the so-called \emph{Inversion Principle} (as put formulated by Gentzen) has the same status as explanation of meaning in the form of movements within language, as do these game/dialogical approaches.\footnote{A number of technical papers and theses have come out of such perspective, as, e.g.\ (Gabbay and de Queiroz 1992), (de Oliveira and de Queiroz 1999, 2005), (de Queiroz 1987,1988ab,1989,1991,1992,1994, 2001,2008,2023,2025), (Mart\'{\i}nez-Rivillas 2022), (Mart\'{\i}nez-Rivillas and de Queiroz 2022ab,2023ab,2025), (de Queiroz and Gabbay 1994,1995,1997,1999), (de Queiroz and Maibaum 1990,1991), (de Queiroz and de Oliveira 2011,2014), (de Queiroz et al.\ 2011,2016), (Ramos 2018), (Ramos et. al 2017,2018,2021ab), (de Veras et al.\ 2021,2023ab).}

\paragraph*{Acknowledgements.} Thanks to the anonymous reviewer(s) of the previous two papers for being extremely careful and sharp enough to contribute decisively to make the text clearer and as more understandable as possible.  

\section*{References}
Church, A.1932.  A set of postulates for the foundation of logic. \textit{Annals
of Mathematics}, Series 2, 33:346--366, 1932.\\
\ \\
Church, A. 1933. A set of postulates for the foundation of logic (second
paper). {\em Annals of Mathematics}, Series 2, 34:839--864, 1933.\\
\ \\
Church, A. 1936. An Unsolvable Problem of Elementary Number Theory. \textit{American Journal of Mathematics}, Vol. 58, No. 2. (Apr., 1936), pp. 345--363, 1936.\\
\ \\
Dummett, M. 1991. \textit{The Logical Basis of Metaphysics}, Harvard University Press, Cambridge (Mass.).\\
\ \\
Frege, G. 1893. \emph{Grundgesetze der Arithmetik}. Hermann Pohle, Jena 1893 (Band I) 1903 (Band II)\\
\ \\
Frege, G. 1965. {\em Basic Laws of Arithmetic}, partial translation by M. Furth, Univ. California Press, 1965.\\
\ \\
Frege, G. 2013. {\em Basic Laws of Arithmetic}, Ph. A. Ebert and M. Rossberg (eds. and trans.), Oxford University Press, 2013.\\
\ \\
Gabbay, D.M., de Queiroz, R.J.G.B. 1992. ``Extending the Curry-Howard Interpretation to Linear, Relevant and Other Resource Logics" \emph{Journal of Symbolic Logic} 57(4):1319--1365.\\
\ \\
Gentzen, G. 1935. Untersuchungen \"uber das logische Schlie{\ss}en. I. \emph{Math Z} 39, 176--210. https://doi.org/10.1007/BF01201353. (English translation in J.\ von Plato (2017)).\\
\ \\
Hintikka,  J. 1968. ``Language-games for quantifiers". In N.\ Rescher, editor, \emph{Studies in Logical Theory}, pages 46--72. Blackwell.\\
\ \\
Hintikka,  J. 1979. ``Quantifiers vs. Quantification Theory". In E. Saarinen (ed) \textit{Game-Theoretical Semantics}. Synthese language library, vol 5. Springer, Dordrecht.\\
\ \\
Hintikka, M. and Hintikka, J. 1986. \textit{Investigating Wittgenstein}.  Basil Blackwell, Oxford.\\
\ \\
Hodges, W. 2001. ``A Sceptical Look". \textit{Proceedings of the Aristotelian Society} Supplementary Volume {75}(1):17--32.\\
\ \\
Howard, W. 1980. ``The formulae-as-types notion of construction".
In H. Curry, J.R. Hindley, J. Seldin (eds.), \emph{To H. B. Curry: Essays on Combinatory Logic, Lambda Calculus, and Formalism}. Academic Press (1980). (Original manuscript circulated in 1969.)\\
\ \\
Lorenzen, P. 1950. ``Konstruktive Begr\"{u}ndung der Mathematik". \emph{Mathematische Zeitschrift} 53(2):162--202.\\
\ \\
Lorenzen, P. 1955. \textit{Einf\"{u}hrung in die operative Logik und Mathematik}. Die Grundlehren der mathematischen Wissenschaften, vol. 78. Springer-Verlag, Berlin-G\"ottingen-Heidelberg. VII + 298 pp\\
\ \\
Lorenzen, P. 1969. \textit{Normative Logic and Ethics\/}, series B.I-Hochschultaschenb\"ucher. Systematische Philosophie, vol. 236$^*$, Bibliographisches Institut, Mannheim/Z\"urich.\\
\ \\
Martin-L\"of, P. 1984. \emph{Intuitionistic Type Theory}. (Notes taken by G.\ Sambin). Bibliopolis, Napoli. ISBN 978-8870881059.\\
\ \\
Martin-L\"of, P. 1987. ``Truth of a Proposition, Evidence of a Judgement, Validity of a Proof". \textit{Synthese} {73}:407--420.\\
\ \\
Martin-L\"of, P.  2013. ``Verificationism Then and Now". Chapter 1 of M.\ van der Schaar (ed.), \emph{Judgement and the Epistemic Foundation of Logic},  Logic, Epistemology, and the Unity of Science 31, DOI 10.1007/978-94-007-5137-8\_1, Springer 2013.\\
\ \\
Martin-L\"of, P.  2019. ``Logic and Ethics". In Piecha, T.; Schroeder-Heister, P. (ed.), \emph{Proof-Theoretic Semantics: Assessment and Future Perspectives}. Proceedings of the Third T\"ubingen Conference on Proof-Theoretic Semantics, 27--30 March 2019,  Univ T\"ubingen, pp.\ 227--235. http://dx.doi.org/10.15496/publikation-35319.\\
\ \\
Martin-L\"of, P.  2022. ``Corrections of Assertion and Validity of
Inference".  Transcript of a lecture given on 26 October 2022 at the Rolf Schock Symposium in Stockholm. \\https://pml.flu.cas.cz/uploads/PML-Stockholm26Oct22.pdf\\
\ \\
Mart\'{\i}nez-Rivillas, D.O. 2022. \emph{Towards a homotopy domain theory}. PhD thesis,
CIn-UFPE (November 2022). Centro de Inform\'atica, Universidade Federal
de Pernambuco, Recife, Brazil. https://repositorio.ufpe.br/handle/123456789/49221\\
\ \\
Mart\'{\i}nez-Rivillas, D.O., de Queiroz, R.J.G.B. 2022a. ``$\infty$-Groupoid Generated by an Arbitrary Topological $\lambda$-Model". \emph{Logic J.\ of the IGPL} 30(3):465--488. https://doi.org/10.1093/jigpal/jzab015 . Also arXiv:1906.05729\\
\ \\
Mart\'{\i}nez-Rivillas, D.O., de Queiroz, R.J.G.B. 2022b. ``Towards a Homotopy Domain Theory". \emph{Archive for Mathematical Logic} 62:559--579. Nov 2022. https://doi.org/10.1007/s00153-022-00856-0\\
\ \\
Mart\'{\i}nez-Rivillas, D.O., de Queiroz, R.J.G.B. 2023a. ``The Theory of an Arbitrary Higher $\lambda$-Model". \emph{Bulletin of the Section of Logic}, 52(1), 39--58. https://doi.org/10.18778/0138-0680.2023.11 . Also arXiv:2111.07092\\
\ \\
Mart\'{\i}nez-Rivillas, D.O., de Queiroz, R.J.G.B. 2023b. ``Solving homotopy domain equations". (Submitted for publication.) arXiv:2104.01195\\
\ \\
Mart\'{\i}nez-Rivillas, D.O., de Queiroz, R.J.G.B. 2025. ``The $K_\infty$ Homotopy $\lambda$-Model". (Submitted for publication.) arXiv:2505.07103\\
\ \\
de Oliveira, A.G., de Queiroz, R.J.G.B. 1999. ``A Normalization Procedure for the Equational Fragment of
Labelled Natural Deduction". \emph{Logic J.\ of the IGPL} 7(2):173--215.\\
\ \\
de Oliveira, A.G., de Queiroz, R.J.G.B. 2005. ``A new basic set of proof transformations". In S.\ Artemov, H.\ Barringer, A.\ Garcez, L.\ Lamb, \& J.\ Woods (Eds.), \emph{We will show them! Essays in Honour of Dov Gabbay} (Vol. 2, pp.\ 499--528). London: College Publications.\\
\ \\
Peirce, C. S. 1932. \emph{Collected Papers of Charles Sanders Peirce, Volumes I and II: Principles of Philosophy and Elements of Logic}. C. Hartshorne and
P.\ Weiss (ed.). Harvard Univ Press.\\
\ \\
von Plato, J. 2017. \emph{Saved from the Cellar.
Gerhard Gentzen's Shorthand Notes on Logic and Foundations of Mathematics.} Springer.\\
\ \\
Plotkin, G. 1975. Call-by-name, call-by-value and the $\lambda$-calculus. \emph{Theoretical Computer Science}
Volume 1, Issue 2, December 1975, Pages 125--159\\
\ \\
Prawitz, D. 1965. \textit{Natural Deduction: A Proof-Theoretical Study}. Acta Universitatis Stockholmiensis, Stockholm Studies in Philosophy no.\ 3. Almqvist \& Wiksell, Stockholm, G\"oteborg, and Uppsala.\\
\ \\
Prawitz, D. 1977. ``Meaning and proofs: on the conflict between classical and intuitionistic logic". \textit{Theoria} (Sweden) {XLIII} 2--40.\\
\ \\
Prawitz, D.  2019. ``Validity of Inferences Reconsidered". In Piecha, T.; Schroeder-Heister, P. (ed.), \emph{Proof-Theoretic Semantics: Assessment and Future Perspectives}. Proceedings of the Third T\"ubingen Conference on Proof-Theoretic Semantics, 27--30 March 2019,  Univ T\"ubingen, pp.\ 213--226. http://dx.doi.org/10.15496/publikation-35319.\\
\ \\
de Queiroz, R.J.G.B. 1987. ``Note on Frege's notions of definition and the relationship proof theory vs.\ recursion theory (Extended Abstract)". In \emph{Abstracts of the VIIIth International Congress of Logic, Methodology and Philosophy of Science.} Vol.\ 5, Part I, Institute of Philosophy of the Academy of Sciences of the USSR, Moscow, 1987, pp.\ 69--73.\\
\ \\
de Queiroz, R.J.G.B. 1988a. ``A Proof-Theoretic Account of Programming and the Role of Reduction Rules". \emph{Dialectica} 42(4):265--282.\\
\ \\
de Queiroz, R.J.G.B. 1988b. ``The mathematical language and its semantics: to show the consequences of a proposition is to give its meaning". In P.\ Weingartner, G.\ Schurz, E.\ Leinfellner, R.\ Haller, A.\ H\"ubner, (eds.) \emph{Reports of the thirteenth international Wittgenstein symposium} 18, pp.\ 259--266.\\
\ \\
de Queiroz, R.J.G.B. 1989. ``Meaning, function, purpose, usefulness, consequences -- interconnected concepts (abstract)". In \emph{Abstracts of Fourteenth International Wittgenstein Symposium (Centenary Celebration)}, 1989, p.\ 20. Symposium held in Kirchberg/ Wechsel, August 13--20 1989.\\
\ \\
de Queiroz, R.J.G.B. 1991. ``Meaning as grammar plus consequences". \emph{Dialectica} 45(1):83--86.\\
\ \\
de Queiroz, R.J.G.B. 1992. ``Grundgesetze alongside Begriffsschrift (abstract)", in \emph{Abstracts of Fifteenth International Wittgenstein Symposium}, 1992, pp.\ 15--16. Symposium held in Kirchberg/Wechsel, August 16--23 1992.\\
\ \\
de Queiroz, R.J.G.B. 1994. ``Normalisation and Language Games". \emph{Dialectica} 48(2):83--123.\\
\ \\
de Queiroz, R.J.G.B. 2001. ``Meaning, function, purpose, usefulness, consequences -- interconnected concepts". \emph{Logic J of the IGPL} 9(5):693--734.\\
\ \\
de Queiroz, R.J.G.B. 2008. ``On Reduction Rules, Meaning-as-use, and Proof-theoretic Semantics". \emph{Studia Logica} 90:211--247.\\
\ \\
de Queiroz, R. J. G. B. 2023.``From Tractatus to Later Writings and Back -- New Implications from Wittgenstein's Nachlass" SATS, vol. 24, no.\ 2,  pp.\ 167--203. https://doi.org/10.1515/sats-2022-0016\\
\ \\
de Queiroz, R. J. G. B. 2025.``From the Notebooks to the Investigations and Beyond" SATS, (to appear). Preprint arXiv:2504.18949\\
\ \\
de Queiroz, R.J.G.B., Gabbay, D.M. 1994. ``Equality in Labelled Deductive Systems and the functional interpretation of propositional equality". In P.\ Dekker, and M.\ Stokhof, (eds.), \emph{Proceedings of the 9th Amsterdam Colloquium 1994}, ILLC/Department of Philosophy, University of Amsterdam, pp.\ 547--546.\\
\ \\
de Queiroz, R.J.G.B., Gabbay, D.M.  1995. ``The functional interpretation of the existential quantifier"'. \emph{Bull of the IGPL} 3(2--3):243--290.\\
\ \\
de Queiroz, R.J.G.B., Gabbay, D.M.  1997. ``The functional interpretation of modal necessity". In M.\ de Rijke, (ed.), \emph{Advances in Intensional Logic}, Applied Logic Series, Kluwer, 1997, pp.\ 61--91.\\
\ \\
de Queiroz, R.J.G.B., Gabbay, D.M. 1999. Labelled Natural Deduction. In: Ohlbach, H.J., Reyle, U. (eds) \emph{Logic, Language and Reasoning}. Trends in Logic, vol 5. Springer, Dordrecht. https://doi.org/10.1007/978-94-011-4574-9\_10\\
\ \\
de Queiroz, R.J.G.B., Maibaum, T.S.E. 1990. ``Proof Theory and Computer Programming". \emph{Zeitschrift f\"ur mathematische Logik und Grundlagen der Mathematik} 36:389--414. \\
\ \\
de Queiroz, R.J.G.B., Maibaum, T.S.E. 1991. ``Abstract Data Types and Type Theory: Theories as Types". \emph{Zeitschrift f\"ur mathematische Logik und Grundlagen der Mathematik} 37:149--166. \\
\ \\
de Queiroz, R.J.G.B., de Oliveira, A.G. 2011. ``The Functional Interpretation of Direct Computations". \emph{Electronic Notes in Theoretical Computer Science} 269:19--40.\\
\ \\
de Queiroz, R.J.G.B., de Oliveira, A.G. 2014.  ``Natural Deduction for Equality: The Missing Entity". In Pereira, L., Haeusler, E., de Paiva, V. (eds) \emph{Advances in Natural Deduction}. Pages 63--91. Trends in Logic, vol 39. Springer, Dordrecht..\\
\ \\
de Queiroz, R.J.G.B., de Oliveira, A.G., Gabbay, D.M. 2011.  \emph{The Functional Interpretation of Logical Deduction}. Vol. 5 of Advances in Logic series. Imperial College Press / World Scientific, Oct 2011.\\
\ \\
de Queiroz, R.J.G.B., de Oliveira, A.G., Ramos, A.F. 2016. ``Propositional Equality, Identity Types, and Computational Paths". \emph{South Amer. J. of Logic} 2(2):245--296.\\
\ \\
Ramos, A.F. 2018. \emph{Explicit computational paths in type theory}. PhD thesis,
CIn-UFPE (August 2018). Centro de Inform\'atica, Universidade Federal
de Pernambuco, Recife, Brazil. \linebreak https://repositorio.ufpe.br/handle/123456789/32902 (Abstract in: Ramos, A. (2019). Explicit Computational Paths in Type Theory. \emph{Bulletin of Symbolic Logic}, 25(2):213-214. \linebreak doi:10.1017/bsl.2019.2)\\
\ \\
Ramos, A.F., de Queiroz, R.J.G.B., de Oliveira, A.G.  2017. ``On the identity type as the type of computational paths". \emph{Logic J.\ of the IGPL} 25(4):562--584.\\
\ \\
Ramos, A.F., de Queiroz, R.J.G.B., de Oliveira, A.G., de Veras, T.M.L. 2018. ``Explicit Computational Paths". \emph{South Amer.\ J.\ of Logic} 4(2):441--484.\\
\ \\
Ramos, A.F., de Queiroz, R.J.G.B., de Oliveira, A.G. 2021a. ``Computational Paths and the Fundamental Groupoid of a Type". In: \emph{Encontro de Teoria da Computa\c{c}\~ao (ETC)}, 6, Evento Online. Porto Alegre: Sociedade Brasileira de Computa\c{c}\~ao, 2021. p.\ 22--25. ISSN 2595-6116. DOI: https://doi.org/10.5753/etc.2021.16371\\
\ \\
Ramos, A.F., de Queiroz, R.J.G.B., de Oliveira, A.G. 2021b. ``Convers\~ao de Termos, Homotopia, e Estrutura de Grup\'oide". In: \emph{Workshop Brasileiro de L\'ogica (WBL)}, 2, Evento Online. Porto Alegre: Sociedade Brasileira de Computa\c{c}\~ao, 2021. p.\ 33--40. ISSN 2763-8731. DOI: https://doi.org/10.5753/wbl.2021.15776\\
\ \\
Scott, D. 1970. Lambda Calculus: some models, some philosophy. In {\em The Kleene Symposium}, J.\ Barwise, H.\ J.\ Keisler \& K.\ Kunen (eds.), North-Holland, 1980, pp. 223--265.\\
\ \\
V\"an\"a\"anen, J. 2022. {The Strategic Balance of Games in Logic}. arXiv:2212.01658\\
\ \\
Wittgenstein, L. 1974. \textit{Letters to Russell, Keynes and Moore},
Ed. with an Introd. by G.\ H.\ von Wright, (assisted by B.\ F.\ McGuinness), Basil Blackwell, Oxford.\\
\ \\
Wittgenstein, L. 1961. \textit{Notebooks 1914--1916,} G.\ H.\ von Wright and G.\ E.\ M.\ Anscombe (eds.), Oxford: Blackwell.\\
\ \\
Wittgenstein, L. 1953. \textit{Philosophical Investigations,} G.E.M.\ Anscombe and R.\ Rhees (eds.), G.E.M.\ Anscombe (trans.), Oxford: Blackwell.\\
\ \\
Wittgenstein, L. 1971. \textit{ProtoTractatus--An Early Version of Tractatus Logico- Philosophicus}, B.\ F.\ McGuinness, T.\ Nyberg, G.\ H.\ von Wright (eds.), D.\ F.\ Pears and B.\ F.\ McGuinness (trans.), Ithaca: Cornell University Press.\\
\ \\
Wittgenstein, L. 1922. 
\emph{Tractatus Logico-Philosophicus},  Translated by C.K.\ Ogden.
Kegan Paul, Trench, Trubner \& Co., Ltd. New York: Harcourt, Brace \& Company, Inc.\\
\ \\
Wittgenstein, L. 2016. Wittgenstein, Ludwig: Interactive Dynamic Presentation (IDP) of Ludwig Wittgenstein's philosophical \emph{Nachlass} [wittgensteinonline.no]. Edited by the Wittgenstein Archives at the University of Bergen (WAB) under the direction of Alois Pichler. Bergen: Wittgenstein Archives at the University of Bergen 2016-. http://wab.uib.no/ (accessed 04 October 2024)

\end{document}